\documentclass[final,oneeqnum,onethmnum,onefignum,onetabnum]{siamltex1213}

\usepackage{amsmath}
\usepackage{amssymb}
\usepackage{booktabs}
\usepackage{thmtools}
\usepackage{tikz}
\usepackage[section]{placeins}

\usepackage{subcaption}
\usepackage{hyperref}
\usepackage{enumitem}
\usepackage{wrapfig}
\usepackage[english, american]{babel}
\usepackage[textsize=tiny]{todonotes}
\usepackage{comment}
\usetikzlibrary{patterns}
\usepackage{cite}
\usepackage{geometry}

\newcommand{\term}[1]{\textbf{#1}}
\newcommand{\tech}[1]{\textsf{#1}}
\newcommand{\op}[1]{\ensuremath{\operatorname{#1}}}

\newcommand{\boldq}{\mathbf{q}}
\newcommand{\bolda}{\mathbf{a}}
\newcommand{\bolde}{\mathbf{e}}

\newcommand{\boldx}{\mathbf{x}}
\newcommand{\boldu}{\mathbf{u}}
\newcommand{\boldv}{\mathbf{v}}

\newcommand{\boldz}{\mathbf{z}}
\newcommand{\boldzero}{\boldsymbol{0}}
\newcommand{\Z}{\mathbb{Z}}
\newcommand{\R}{\mathbb{R}}
\newcommand{\C}{\mathbb{C}}
\newcommand{\Cstar}{(\mathbb{C} \setminus \{0\})}

\DeclareMathOperator{\nvol}{NVol}

\DeclareMathOperator{\conv}{conv}

\newtheorem{remark}{Remark}
\newtheorem{example}{Example}
\newtheorem{problem}{Problem Statement}

\begin{document}
\begin{flushright}
 Preprint: ADP-15-50/T952
\end{flushright}

\title{On the Network Topology Dependent Solution Count of the Algebraic Load Flow Equations}

\author{
    Tianran Chen\footnotemark[2] \and
    Dhagash Mehta\footnotemark[3]\ \footnotemark[4]
}

\maketitle

\renewcommand{\thefootnote}{\fnsymbol{footnote}}
\footnotetext[2]{Dept.~of Mathematics, Michigan State University, East Lansing, MI USA.}
\footnotetext[3]{Dept.~of Applied and Computational Mathematics and Statistics, University of Notre Dame, Notre Dame, IN 46556, USA.}
\footnotetext[4]{Centre for the Subatomic Structure of Matter, Department of Physics, School of Physical Sciences, University of Adelaide, Adelaide, South Australia 5005, Australia.}
\renewcommand{\thefootnote}{\arabic{footnote}}

\begin{abstract}
	\noindent
	%Power flow computation is an important problem in the design and control 
	%of power networks.
	A large amount of research activity in power systems areas has focused 
    on developing computational methods to solve load flow equations
	where a key question is the maximum number of isolated solutions.
	Though several concrete upper bounds exist, recent studies have hinted 
    that much sharper upper bounds that depend the topology of underlying 
    power networks may exist.
	This paper establishes such a topology dependent solution bound which is 
    actually the best possible bound in the sense that it is always attainable. 
    We also develop a geometric construction called adjacency polytope which 
    accurately captures the topology of the underlying power network and is
    immensely useful in the computation of the solution bound.
    Finally we highlight the significant implications of the development of
    such solution bound in solving load flow equations.
    %We also simulataneously develop an associated AP 
    %numerical polynomial homotopy continuation method that guarantees to find all isolated complex solutions, and 
    %provide results for plenty of networks.
\end{abstract} 
%\linenumbers
%===============================================================================
\section{Introduction} \label{sec:intro}

Engineers are regularly required to perform power flow computations for designing, 
operating, and controlling power systems \cite{kundur1994power}. 
The solutions of the power flow or \textit{load flow equations}, a system of 
\textit{multivariate nonlinear} equations, are used to ensure the normal 
operating conditions as well as to perform contingency, stability, and 
bifurcations analysis. 
In addition, these solutions, which we refer to as 
\textit{solutions of the load flow equations} from now on, 
are also important in system security assessment and optimal dispatching.
In general, load flow equations may have more than one solutions \cite{klos1975}.
There are indeed quite a few existing methods for finding one or many 
load flow solutions \cite{iwamoto1981,glover_sarma_overbye,thorp1989,
bromberg2015,trias2012,ilic1992,dorfler2013,dj2015,lavaei_convexpf,
dreesen2009,lasserre_book,ajjarapu1992continuation,thorp1993,chen2011,counterexample2013,liu2005,
mori1999,montes1995solving,ning2009,nguyen2014appearance,parrilo2014,sommese_numerical_2005,li_numerical_2003,mehta2014a,
mehta2014b,chandra2015equilibria,mehta2015algebraic} (see \cite{mehta2015recent} for a recent review). 
Out of the few methods that guarantee to find \emph{all} load flow solutions, 
i.e., the interval based approach \cite{mori1999}, 
Gr\"obner bases technique \cite{montes1995solving,ning2009,nguyen2014appearance,parrilo2014}
and the numerical polynomial homotopy continuation (NPHC) method
\cite{sommese_numerical_2005,li_numerical_2003,mehta2014a,mehta2014b,chandra2015equilibria,mehta2015algebraic},
the NPHC method appears most promising in scalability with increasing system 
sizes in that it has already found all load flow solutions of up to IEEE 14 
bus systems \cite{mehta2014a} 
(and 18 oscillators case for the Kuramoto model \cite{mehta2015algebraic}) 
and is inherently parallel: formulating load flow equations as system of 
polynomial equations, the NPHC method, whose roots are in complex algebraic geometry, finds 
all isolated complex systems which obviously include all the isolated real solutions.
In all these computational methods, the knowledge of the number of solutions
play a crucially important role.
%Here, to overcome the scalability issue of this method, one needs to 
%come up with a tight upper bound on the number of isolated complex solutions. 
%Apart from the NPHC method, if such a tight upper bound is known, then one can 
%ensure to have found all solutions using any 
%other numerical method as the upper bound are found if as exactly many solutions.

While several concrete upper bounds exist, several studies \cite{guo_determining_1990,molzahn-mehta-matt2015} has hinted 
that much sharper upper bounds that depend the topology of underlying 
power networks may exist.
The main goal of this paper is to establish such a topology dependent solution bound. 
The remainder of the paper is organized as follows: 
In Sec.~\ref{sec:algebraic}, we set up a formulation of the algebraic load flow 
equations, define the problems precisely, and provide a brief historical account 
on the available results.
Sec.~\ref{sec:Max_No_Sols} describes our rigorous results: 
in Sec.~\ref{sec:bernshtein} we describe the tight bound, called the BKK bound, 
on the number of isolated complex solutions for our systems. 
We also refine the BKK bound in the context of algebraic load flow systems.
%In particular, we show the BKK bound is always attainable even for the quite special
%load flow systems where
In Sec.~\ref{sec:polytope-bound} we introduce a novel upper bound called 
adjacency polytope (AP) bound. 
Sec.~\ref{sec:computing-bound} discusses the computational issues,
and in Sec.~\ref{sec:sol_count_homotopy} we highlight the significant 
implications of the development of these solution bounds in homotopy methods 
for solving load flow equations.
In Sec.~\ref{sec:cases}, we compute the BKK and AP bounds and the actual number 
of complex solutions for numerous power network examples from literature and 
compare with the previously known bounds. 
In Sec.~\ref{sec:conclusion}, we discussion our results 
in a larger context and conclude.

%===============================================================================
\section{Algebraic formulation}\label{sec:algebraic}

In this paper, we focus on the mathematical abstraction of a power network which 
is captured by a graph $G = (B,E)$ together with a complex matrix $Y = (Y_{ij})$.
Here $B$ is the finite set of nodes representing the ``buses'', $E$ is the 
set of edges (also called lines or branches) representing the connections between buses, and the matrix 
$Y$ is the nodal \textit{admittance matrix}\footnote{%
	As such, the nodal admittance matrix $Y$ already captures the graph topology
	since its sparsity pattern describes precisely the edges and hence the 
    topology of the graph $G=(B,E)$.
    However, since we often make separate use of the graph topology
    and the nodal admittance information, it is therefore convenient to 
    explicitly extract the graph structure $G$ and keep it separate from $Y$.
} which assigns a nonzero complex value $Y_{ij}$
(mutual admittances) to each edge in $(i,j) \in E$.
For any $(i,j) \notin E$ (i.e., nodes $i$ and $j$ are not directly connected),
$Y_{ij} = Y_{ji} = 0$. 
Here, $Y$ is \emph{not} assumed to be symmetric, but we require $Y_{ij}$ and 
$Y_{ji}$ to be either both zero or both nonzero (i.e., the underlying graph 
is ``undirected'').
As a convention, we further require all nodes to be connected with itself via a
``loop'' to reflect the nonzero diagonal entries $Y_{ii}$ known as the 
\emph{self-admittances}.
That is, we require $(i,i) \in E$ for any $i \in B$.

For brevity, we define $n$ to be the number of non-reference buses
(i.e., $|B| = n + 1$) and label the nodes by integers $0,1,2,\dots,n$. 
The corresponding complex valued voltages will be denoted by $v_0, v_1,\dots,v_n$.
Here we fix node 0 to be the designated reference bus for which the voltage is 
fixed to a nonzero real number, that is, $v_0 \in \R$ is a nonzero constant.
%With this notation we consider the (weighted) adjacency matrix $Y = (Y_{jk})$ 
%where each entry $Y_{jk}$ is the complex weight of the edge between nodes 
%$i$ and $j$ if such an edge exists and $0$ otherwise.
%We do not require this matrix to be symmetric.
In this setup, the load flow equation takes the form of
\begin{equation}
    \label{equ:Sj}
    S_i = \sum_{j=0}^n Y_{ij}^* v_i v_j^* \quad \text{for } i = 1,\dots,n,
\end{equation}
which is a system of $n$ equations in the $n$ variables $v_1,\dots,v_n$ since 
$v_0$, corresponding to the reference bus, is fixed (and hence not a variable).
Here $v_i^*$ and $Y_{ij}^*$ denotes the complex conjugates of $v_i$ and $Y_{ij}$
respectively, and $S_i\in \mathbb{C}$ are the injected power. 
The equations \ref{equ:Sj} may represent either a transmission or a distribution network,
with PQ buses.
It is the network topology along with other features that can distinguish these cases: 
a mesh topology would usually correspond to transmission networks, whereas radial (tree-like) 
topology would correspond to distribution networks.

A solution to \eqref{equ:Sj} is said to be \term{isolated} if it is the only 
solution in a sufficiently small neighborhood.
i.e., an isolated solution has no degree of freedom.
Solutions with some $v_k = 0$ are said to be \term{deficient}, and 
\term{non-deficient} otherwise.
By a standard application of Sard's Theorem\cite{sommese_numerical_2005}, it can be verified that under a
generic perturbation of $S_1,\dots,S_n$, the system \eqref{equ:Sj} has
no deficient solutions (though not a completely unlikely case in some load flow systems).
We therefore focus only on the non-deficient solutions.
%For a deficient solution with a some $v_k = 0$, one can remove all terms 
%involving $v_k$ and $v^*_k$ from \eqref{equ:Sj} and equations would still 
%balance as if the $k$-th bus are absent. 
%In other words, a deficient solution represents a solution to the subnetwork with the 
%$k$-th bus removed.
The problem central to this paper is a counting problem for the isolated 
non-deficient load flow solutions:
%and we ignore the physically less likely deficient solutions\footnote{It may however
%be interesting to conduct a similar study in which deficient solutions are included in future.}:

\medskip
\begin{problem}
    \label{prb:couting}
    For a power network, what is the maximum number of isolated 
    non-deficient solutions to the system \eqref{equ:Sj}?
\end{problem}
\medskip

Following the fruitful algebraic approach taken by works such as 
\cite{baillieul1982,li_numerical_1987}, we shall ``embed'' this problem 
into a bigger algebraic root counting problem: 
That is, we consider a polynomial system whose solution set captures 
all the solutions of the above (non-algebraic) system.
To that end, we introduce new variables
\begin{equation}
    \label{equ:conjugate}
    u_k = v_k^* \quad \text{for each } k = 0, \dots, n.
\end{equation}
Substituting them into \eqref{equ:Sj}, we obtain the system of algebraic equations
\begin{equation}
    \label{equ:Sj-alg}
    S_j = \sum_{k = 0}^n Y_{jk}^* \, v_j \, u_k 
    \quad \text{for } j = 1,\dots,n.
\end{equation}
This is a system of $n$ equations in $2n$ variables.
However, a ``square'' system where the number of variables and equations match
is often much more convenient from an algebraic point of view.
We therefore extract $n$ hidden equations by taking the complex conjugate of 
both sides of each of the above equations and obtain
\begin{equation}
    \label{equ:Sj-conj}
    S_j^* = 
    \sum_{k=0}^n (Y_{jk}^* \, v_j \, u_k)^* =
    \sum_{k=0}^n  Y_{jk}   \, u_j \, v_k 
         \quad \text{for } j =1,\dots,n.
\end{equation}
%since we require $u_k = v_k^*$.
We now sever the tie between $\boldu=(u_1,\dots,u_n)$ and $\boldv=(v_1,\dots,v_n)$ 
and consider them to be variables independent from one another.
That is, we ignore \eqref{equ:conjugate}.
Then \eqref{equ:Sj-alg} and \eqref{equ:Sj-conj} combine into a system of
$2n$ polynomial equations in the $2n$ variables: 
%$v_1,\dots,v_n,u_1,\dots,u_n$
%which can be written as
\begin{equation}
    P_{G,Y} (v_1,\dots,v_n,u_1,\dots,u_n) = 
    \begin{cases}
        \sum_{k=0}^{n} Y_{1k}^*\, v_1 u_k - S_1   = 0\\
        \quad \vdots \\
        \sum_{k=0}^{n} Y_{nk}^*   v_n u_k - S_n   = 0\\[1.5ex]
        \sum_{k=0}^{n} Y_{1k}\,\, u_1 v_k - S_1^* = 0 \\
        \quad \vdots \\
        \sum_{k=0}^{n} Y_{nk}\,   u_n v_k - S_n^* = 0.
    \end{cases}
    \label{equ:powerflow-alg}
\end{equation}
Note that here the values of $v_0$ and $u_0$ are fixed, as they correspond 
to the reference node and are hence constants in the above system.
For brevity, this system will be referred to as the 
\term{algebraic load flow equations} in the following discussions.
This formulation first appeared in \cite{li_numerical_1987} to the best of our knowledge.
Other polynomial formulations of the load flow equations have been also 
known (see, e.g., \cite{baillieul1982,baillieul1984critical,mehta2014a,mehta2014b,chandra2015equilibria,marecek2014power}).
    
It is worth noting that in $P_{G,Y}$, the topology (i.e., the edges) of the 
underlying power network, in a sense, is encoded in the set of monomials 
while entries of $Y$ appear as the coefficients.
Developing a solution count that exploits network topology via the monomial
structure is the main goal of the current article.
%and this system has sparked important development in the community of
%\emph{homotopy continuation method}.
    
%\begin{remark}
%    The formulation \eqref{equ:powerflow-alg} first appeared in 
%    \cite{li_numerical_1987}.
%    A similar polynomial formulation for a special case was also proposed in
%    earlier works \cite{baillieul1982,baillieul1984critical}.
%    It is, however, worth noting that this formulation is a the deeper motivation behind this transformation.
%    A basic observation is that in a function $f(z)$ in the complex variable $z$
%    (which may or may not be holomorphic), the variable $z$ and its conjugate 
%    $z^*$ can be treated as independent variables.
%    That is, in certain context, one can treat $f$ as the function $f(z,z^*)$.
%    In the German-speaking world, this technique is commonly known as ???.
%    \cite{???} provides a thorough study of this technique.
%    The formulation \eqref{equ:powerflow-alg} is therefore the load flow equation
%    in the $z$-$z^*$ coordinate system.
%\end{remark}

Clearly, for every solution $\boldv = (v_1,\dots,v_n)$ of the original 
(non-algebraic) system \eqref{equ:Sj}, %by construction,
we naturally have $P_{G,Y}(\boldv,\boldv^*) = \boldzero$.
That is, the solution set of $P_{G,Y} = \boldzero$ captures all solutions of 
the original (nonalgebraic) load flow system \eqref{equ:Sj}.
In the following we shall focus on the following root counting problem:

\medskip
\begin{problem}
    \label{prb:alg-counting}
    Among the power networks of a given topology provided by a graph $G$, 
    what is the maximum number of isolated solutions of $P_{G,Y}$ in $\Cstar^{2n}$
    for all possible choices of $Y$?
\end{problem}
\medskip

Here, the phrase ``maximum number'' means the 
\term{lowest upper bound that is also attainable} 
and shall be distinguished from a mere ``upper bound''.
Of course, the existence of such a ``maximum number'' is not a priorily guaranteed,
after all the lowest upper bound taken over an infinite family (choices of $Y$'s),
even if exists, may not be attainable.
One of the goals of this paper is to establish the validity of the above question.
Clearly, any answer to Problem \ref{prb:alg-counting} provides an upper bound
for the answer to Problem \ref{prb:couting}.

Various upper bounds for Problem \ref{prb:alg-counting} have been proposed 
in the past (see \cite{molzahn-mehta-matt2015} for a recent review). 
The \emph{classical B\'ezout number} (CB number) provides a simple
upper bound: it is the product of the degrees of each polynomial equation in 
\eqref{equ:powerflow-alg} (hence also known as the \emph{total degree}).
It is a basic fact in algebraic geometry that the number of isolated complex
solutions of a polynomial system is bounded above by the CB number.
Therefore, for a power network of $n$ (non-reference) buses, and one reference bus, 
the CB bound is $2^{2n}$, since there are $2n$ equations in \eqref{equ:powerflow-alg}
each of degree 2.
A much tighter upper bound on the number of isolated complex solutions, 
$\binom{2n}{n}$, was derived for the special case of completely
interconnected lossless networks with all nodes being power-voltage nodes in Ref.~\cite{baillieul1982} and 
then for the general case in Ref.~\cite{li_numerical_1987} (see \cite{marecek2014power} for a recent alternative derivation
of this bound). We shall refer to this bound 
as the Baillieul-Byrne-Li-Sauer-Yorke 
(BBLSY) bound.
However, neither of these bounds exploit the network topology of a given power system.
The link between network topology and complex solution count was first hinted in
\cite{guo_determining_1990}, however, a concrete and computable answer remains elusive.
%The network topology of a specific network, networks with cliques with exactly one common node,
%was exploited in \cite{guo1990} and the number (not just an upper bound) of complex solutions was shown be equal 
%to the product of number of complex solutions for each clique as independent network. Recently, 
%the result was extended to networks with cliques sharing exactly two common nodes in \cite{molzahn-mehta-matt2015}.
%For real solutions, there is rarely a precise work available except 
%\cite{chiang2015} which computed the maximum number of type-1
%unstable solutions to be $2\binom{2n_g}{n_g}$ with $n_g$ being the number of generator buses, though 
%under certain conditions. 

In a recent study \cite{molzahn-mehta-matt2015}, 
with extensive numerical experiments via the NPHC methods, 
it was observed that the number of isolated complex solutions is generally 
significantly lower than both the CB and BBLSY bound for sparsely connected graphs. 
Based on these observations, it was anticipated that the key to exploiting the 
network structure of the power system may be to exploit the underlying topology
of the power system. 
%It is this observation that we exploit here and produce a systematic way to 
%compute the maximum number of isolated complex solutions for any given
%network structure.
In the present work, we show that this maximum number exists and it is essentially 
the \emph{Bernshtein-Kushnirenko-Khovanskii} (or BKK) bound.
%\cite{bernshtein_number_1975,kushnirenko_newton_1976,khovanskii_newton_1978}. 
%As a theoretical exploration of this problem, we also highlight the rich 
%algebraic structure of the underlying equation. 
We then develop a novel approximation of this maximum number,
the ``adjacency polytope'' (or AP) bound, which has tremendous computational 
advantage over the BKK bound.
Moreover, we show, via experiments, that the AP bound is exactly the BKK bound
for many concrete cases, making it an ideal surrogate for the BKK bound.
%For the algebraic load flow equations, it turns out that the AP is always equal
%to the BKK bound, for the given system with generic parameters, which also turn out to be the 
%actual number of isolated complex solutions.
%In general, it is possible for the system \eqref{equ:powerflow-alg} to have 
%solutions that are not isolated.
%However, as we shall show,  

%It is worth noting the inherent symmetry in the algebraic load flow system
%\eqref{equ:powerflow-alg} which is a result of treating $u_k$ and $v_k$ to be
%independent variables:
%
%\begin{lemma}
%    \label{lem:symmetry}
%    If $(\boldv,\boldu) = (v_1,\dots,v_n,u_1,\dots,u_n) \in \Cstar^{2n}$ is a
%    solution of the system \eqref{equ:powerflow-alg}, then 
%    $(\boldu^*,\boldv^*) = (v_1^*,\dots,v_n^*,u_1^*,\dots,u_n^*)$ is also a
%    solution.
%\end{lemma}
%
%\begin{proof}
%    \todo[inline]{Finish the proof here.}
%\end{proof}

%%===============================================================================
%\section{A short review of the technical tools}
%
%Below we provide a concise and informal review of the technical tools
%needed in the main discussions and proofs.
%The detailed descriptions are included in the appendix for completeness.
%
%\begin{description}
%    \item[Mixed volume]
%    \item[BKK bound]
%    \item[Randomization]
%    \item[Polarization]
%\end{description}

%===============================================================================
\section{The Maximum Number of Solutions}\label{sec:Max_No_Sols}

We shall focus on Problem \ref{prb:alg-counting} %and consider the maximum 
%number of isolated non-deficient solutions for \eqref{equ:powerflow-alg}.
%Various upper bounds have been developed in the past as briefly reviewed in 
%\S \ref{sec:intro}.
and show that the BKK bound
%By applying Bernshtein's theorem \cite{bernshtein_number_1975} to the 
%algebraic load flow equations and demonstrate that it 
provides the best possible bound for any given power network topology 
which is always attainable by some choice of the $Y$-matrix.
Then in \S \ref{sec:polytope-bound} we propose a simplified formulation 
that is easier to compute and analyze.
Computational issues are explored in \S \ref{sec:computing-bound}.

%===============================================================================
\subsection{The BKK bound} \label{sec:bernshtein}

Problem \ref{prb:alg-counting} is a special case of the root counting problem
for systems of polynomial equations which is an important problem in algebraic
geometry and mathematics in general that has wide range of applications \cite{sommese_numerical_2005,morgan_computing_1987,morgan_finding_1989}.
Two basic root counts are provided by the CB and BBLSY bounds described above.
%which considers the ``multi-degree'' derived from certain grouping of the variables.
One common weakness of these two upper bounds is that they only utilize 
the rather incomplete information about the polynomial system --- the degree
(or ``multi-degree'' of the polynomials). 
In the context of the algebraic load flow equation, this means that they do not 
take into consideration the topology of the underlying power network.
Following up from the observations in Ref.~\cite{molzahn-mehta-matt2015}
we refine these bounds using the 
theory of BKK bound \cite{bernshtein_number_1975} which accurately captures 
the network topology of the power systems.
%In other words, this bound takes into consideration the much more refined information: the 
%monomial structure (the monomials that actually appear) of a polynomial system,
%and consequently it provides a much sharper root count that is, in a sense,
%``generically exact''.
Here we state the powerful, albeit abstract, Bernshtein's theorem in the
context of the algebraic load flow equation.

\medskip
\begin{theorem}[Bershtein \cite{bernshtein_number_1975}]
    \label{thm:bernshtein}
    Consider the algebraic load flow system of $2n$ polynomial equations 
    \eqref{equ:powerflow-alg} in $2n$ variables.
    \begin{enumerate}[label=(\Alph*)]
        \item 
            The number of isolated solutions the system has in $\Cstar^{2n}$ is 
            bounded above by the mixed volume of the 
            Newton polytopes of the $2n$ equations.
        \item
            Without enforcing the conjugate relations among the coefficients%
            \footnote{That is, we treat coefficients with ``$*$'' e.g. $Y_{jk}^*$ 
                and $S_k^*$ in \eqref{equ:powerflow-alg} to be independent from 
                their counterparts without ``$*$''.}
            there is a nonempty Zariski open set of coefficients for which all 
            solutions of the system \eqref{equ:powerflow-alg} in $\Cstar^{2n}$ 
            are isolated and the total number is exactly the upper bound given in (A).
    \end{enumerate}
\end{theorem}
\medskip

The technical terms 
%(such as mixed volume, Newton polytopes, and Zariski open sets) 
are explained in Appendix \ref{sec:mvol} for completeness. 
Here, it is sufficient to take the following interpretation:
Part (A) establishes a computable upper bound for the number of isolated solutions
that depends on the geometric configuration of the monomial structure,
and part (B) shows this upper bound is generically exact.%
\footnote{%
    The above theorem states that ignoring the conjugate relations among 
    the coefficients, the set of coefficients for which the BKK bound holds 
    exactly is at least a nonempty ``Zariski open'' set.    
    Such a set is always open and dense and hence of full measure. 
    Consequently, the BKK bound is commonly known as being ``generically exact'': 
    when coefficients are chosen at random, with probability one, 
    the BKK bound is exact.
}
The original proof was given in \cite{bernshtein_number_1975}. 
An alternative proof that gives rise to the development of \emph{polyhedral homotopy}
is developed in \cite{huber_polyhedral_1995}.
More detail can be found in references such as 
\cite{huber_bernsteins_1997,li_numerical_2003,chen_theoretical_2014}.
In \cite{canny_optimal_1991}, the root count in the above theorem was nicknamed 
the \term{BKK bound}\footnote{The BKK bound is also known as the Bernshtein's bound, 
the polyhedral bound, or the mixed volume bound in literature.}
after the works of Bernshtein \cite{bernshtein_number_1975}, 
Kushnirenko \cite{kushnirenko_newton_1976}, and 
Khovanskii \cite{khovanskii_newton_1978}. 
In general, it provides a much tighter bound on the number of isolated zeros 
of a polynomial system compared to variants of B\'ezout bounds.
More importantly, it is sensitive to the monomial structure and hence
the topology of the underlying power network.

It is important to note that the ``generic exactness'' expressed in part (B)
of the above theorem only holds when one ignores the tie between $Y_{ij}$ and 
$Y_{ij}^*$ as well that between $S_i$ and $S_i^*$.
That is, one must allow $Y_{ij}$ and $Y_{ij}^*$ to vary independently in
interpreting the above theorem.
We shall now bring back the restriction that all the $(Y_{ij},Y_{ij}^*)$ and
$(S_i,S_i^*)$ must be conjugate pairs and investigate the sharpness of the 
BKK bound under these restriction.
Since the BKK bound is sensitive to the monomial structure in the equation
\eqref{equ:powerflow-alg} and hence topology of the underlying power network,
we shall fix the network topology in the following discussion.
That is, we shall fix the sparsity pattern of the $Y$ matrix but allow its
entries to vary among the set of nonzero complex numbers. 
In the following, we shall establish that within the set of admittance matrices
of the same sparsity pattern:
\emph{The BKK bound of the algebraic load flow equations is always exact for some choices of the admittance matrix.}
%\begin{enumerate}
%    \item The BKK bound is always exact for some choices of the admittance matrix; and
%    \item Every admittance matrix is ``infinitely close'' to the set of admittance
%        matrices for which the BKK bound is exact.
%\end{enumerate}
%For a more precise disucssion, we let $\mathcal{Y}$ be the set of all possible 
%admittance matrix having the given sparsity pattern, 
%and let $\mathcal{Y}_<$ be those $Y \in \mathcal{Y}$ for which the number of 
%isolated solutions $\eqref{equ:powerflow-alg}$ has in $\Cstar^{2n}$ is strictly 
%less than the BKK bound.
%Then the above two assertions translate to
%\begin{enumerate}
%    \item $\mathcal{Y}_< \ne \mathcal{Y}$; and
%    \item Any choice $Y \in \mathcal{Y}_<$ is infinitely close to the set
%        $\mathcal{Y} \setminus \mathcal{Y}_<$.
%\end{enumerate}
%We summarize these into the theorems below.
%The proofs and related lemmas are included in Appendix \S \ref{sec:proofs}.
In other words, we have the following assertions:

%That is for a $Y \in \mathcal{Y}$, $Y_{ij} \ne 0$ precisely when
%Though in a realistic power network certain restrictions for the entries in the
%$Y$ matrix such as symmetry and limited range are usually in place, here, as a
%These two assertions are proved in Corollaries \ref{cor:attainable} and \ref{cor:generic}
%respectively. Both proofs hinges on the fact that $\mathcal{Y}_<$ is a very
%``sparse'' set that has an empty interior.
%We therefore start with this somewhat abstract characterization of the set 
%$\mathcal{Y}_<$.

%\begin{theorem}
%    Fixing the topology of the underlying power network (i.e., fixing the set
%    of edges and hence the sparsity pattern of $Y$), let $\mathcal{Y}_<$ be the 
%    the set of admittance matrices for which the number of isolated complex 
%    solutions of \eqref{equ:powerflow-alg} in $\Cstar^{2n}$ is strictly less 
%    than the BKK bound. Then $\mathcal{Y}_<$ has no interior point.
%\end{theorem}

\medskip
\begin{theorem}
    \label{thm:attainable}
    Given a graph $G$, there exists an admittance matrix $Y$ on $G$ for which 
    the number of isolated solutions of the corresponding algebraic load flow 
    equation $P_{G,Y}$ in $\Cstar^{2n}$ is exactly the BKK bound.
\end{theorem}
\medskip

\begin{proof}
    For convenience, let $Z = (Z_1,\dots,Z_\ell)$ which collect all the nonzero
    entries of $Y_{jk}$ and $S_j$. 
    That is, $Z$ contains all the nonzero coefficients in \eqref{equ:powerflow-alg}.
    By the interpretation given in Appendix \ref{sec:zariski}, we simply have to
    show that there exists a choice of $Z \in \Cstar^\ell$ such that
    $D(Z,Z^*) \ne 0$.
    
    Suppose no such choice of $Z$ exist, then $D(Z,Z^*) = 0$ for all
    $Z \in \Cstar^\ell$. 
    By Lemma \ref{lem:polarization} in Appendix \ref{sec:polarization}, 
    $D(Z,W) = 0$ for all $(Z,W) \in \Cstar^{2\ell}$. This means $D$ must be 
    a zero polynomial, which is a contradiction.
    Therefore, we must conclude that there is always a choice of $Z$
    (and hence $Y$ and $S$) such that $D(Z,Z^*) \ne 0$.
\end{proof}
\medskip

%\noindent The proof and related lemmas are included in Appendix \S \ref{sec:proofs}.
%\begin{theorem}
%    \label{thm:generic}
%	Given a graph $G = (B,E)$, an admittance matrix $Y$ on $G$, and a threshold 
%    $\epsilon > 0$, there exists an admittance matrix $\tilde{Y}$ on $G$
%    (and hence having the same sparsity pattern as $Y$) such that 
%	$\|Y - \tilde{Y}\| < \epsilon$ and the number of isolated solutions 
%	in $\Cstar^{2n}$ of the algebraic load flow equation $P_{G,\tilde{Y}}$ 
%    in \eqref{equ:powerflow-alg} is exactly the BKK bound.
%\end{theorem}
%\begin{proof}
%	Suppose there exists an open domain $\Omega \subset \C^{n \times |B|}$
%	such that $G(Y,Y^*) = 0$ for all $Y \in \Omega$.
%	Then by the Polarization Lemma (Lemma \ref{lem:polarization}),
%	$G(Y,Z) = 0$ for all $(Y,Z) \in \Omega \times \Omega^*$ where
%	$\Omega = \{ Z \in \C^{n \times |B|} \mid Z^* \in \Omega \}$.
%	That is, $G$ is identically zero on an open domain. A contradiction.
%	Therefore, we must conclude that the set $U$ has no interior point.
%	In other words, every choice of $Y \in U$ is arbitrarily close to some
%	$\tilde{Y} \not\in U$.
%\end{proof}

%===============================================================================
\subsection{Solution bound via adjacency polytope} \label{sec:polytope-bound}

In this section, we develop a significantly simplified version of the BKK bound
for the algebraic load flow equation \eqref{equ:powerflow-alg} that can be
analyzed and computed more easily.
%Indeed, as Remark \ref{rmk:union-benefit} below will discuss, it will also 
%likely to be much easier to compute using numerical methods.
Here we borrow a construction from coding theory and encode the given graph
into a polytope (which, roughly speaking, 
is a geometrical object whose all sides are flat) which we shall call the ``symmetric AP''.
The definition requires the following notations:
We define $\bolde_0 := \boldzero \in \R^n$.
For $i = 1,\dots,n$, $\bolde_i$ denotes the vector in $\R^n$ that has 
an entry 1 on the $i$-th position and zero everywhere else.
$(\bolde_i,\bolde_j)$ is simply the concatenation of the two vectors 
$\bolde_i \in \R^n$ and $\bolde_j \in \R^n$ that forms a vector in $\R^{2n}$.
Finally, ``conv'' denotes the convex hull operator which produces the smallest
convex set containing a given set.

\medskip
\begin{definition}
    \label{def:polytope}
    Given a graph $G=(B,E)$, let
    \begin{equation*}
        \Gamma_G := 
        \bigcup_{(i,j) \in E} \{ (\bolde_i, \bolde_j) \}
        \;\subset \R^{2n} .
    \end{equation*}
    With this, we define the \term{symmetric adjacency polytope} of $G$ to be
%    the network 
%    (represented by the graph $G$) to be the convex hull of $\Gamma_G$ 
%    together with the origin:
    \begin{equation*}
        \nabla_G := \operatorname{conv} (\Gamma_G \cup \{ \boldzero \}).
    \end{equation*}
\end{definition}
%\medskip

Implicit in this definition is the assumption that the graph $G=(B,E)$ is
\emph{undirected}, that is, if $(i,j) \in E$ then $(j,i) \in E$.
The polytope $\nabla_G$ is a geometric encoding of the connectivity of 
the underlying power network with connections manifested as points.
It is only dependent on the connectivity among the nodes
(and hence the sparsity patterns in $Y$)
but independent from the actual entries in $Y$.

\medskip
\begin{remark} %[Computational advantage of $\nabla_G$]
    \label{rmk:union-benefit}
    It is evident in the construction of \eqref{equ:powerflow-alg} that
    equations in the system always contain many common monomials.
    Indeed, if $(i,j)$ is an edge in the graph, then the monomial $v_i u_j$ 
    appear in both the $i$-th and the $(n+j)$-th equation.
    In other words, the monomial structure of \eqref{equ:powerflow-alg} has
    certain level of built-in redundancy.
    Such redundancy is removed in the construction of $\nabla_G$:
    In Definition \ref{def:polytope}, one takes the union of the set of points 
    representing the edges, common monomials in \eqref{equ:powerflow-alg} will 
    therefore coalesce into the same points. 
    Consequently, the polytope $\nabla_G$, in a sense, contains much less
    information than the monomial structure in \eqref{equ:powerflow-alg}.
    Therefore the encoding $\nabla_G$ is advantageous from a computational 
    point of view.
\end{remark}

With the constructions above, we propose a new upper bound on the isolated
non-deficient complex solutions to the algebraic load flow equation that
takes into consideration the connectivity of the underlying power network:

\medskip
\begin{theorem}
    \label{thm:polytope-bound}
    The number of isolated solutions the algebraic load flow system
    \eqref{equ:powerflow-alg} has in $(\C\setminus\{0\})^{2n}$ is 
    bounded above by
    \begin{equation*}
        \mu_G \;:=\; \nvol_{2n} (\nabla_G)
    \end{equation*}
    which we shall define as the the \term{adjacency polytope bound} 
    (or simply \term{AP bound}) for any $P_{G,Y}$ 
    (with any choice of admittance matrix $Y$).
\end{theorem}
\medskip

Here ``$\nvol_{2n}$'' denotes the \emph{normalized volume} in $\R^{2n}$.
It is a volume measurement commonly used in the study of ``lattice polytopes'', 
and it is defined so that the standard ``corner simplex'' (the corner of a
unit hypercube) has volume 1.
\medskip
\begin{proof}
    For a nonsingular $2n \times 2n$ matrix $M$, we can form the new system 
    $M \cdot P_{G,Y}$ by forming the formal matrix-vector product where $P_{G,Y}$ 
    is considered as a column vector.
    This technique is known as \emph{randomization}.
    Clearly, $M \cdot P_{G,Y}(\boldv,\boldu) = \boldzero$ if and only if 
    $P_{G,Y}(\boldv,\boldu) = \boldzero$ and the number of isolated solutions
    (in $\Cstar^{2n}$) remains the same under this transformation.
    It is easy to verify that the \emph{support} of the randomized system
    $M \cdot P_{G,Y}$ is \emph{unmixed} of type $2n$, and the Newton polytope
    is precisely the symmetric adjacency polytope $\nabla_G$ defined in
    Definition \ref{def:polytope}.
    Then by the unmixed form of Bernshtein's Theorem 
     \cite{huber_polyhedral_1995}, the BKK bound of this 
    randomized system is precisely the normalized volume $\nvol_{2n}(\nabla_G)$.
\end{proof}
\medskip

%The proof is included in Appendix \S \ref{sec:proofs}.
%For later references, we shall name the upper bound provided by Theorem 
%\ref{thm:polytope-bound} as ``adjacency polytope bound''.
%The rest of this paper will be devoted to the computation and analysis
%of this bound for various types of power networks.

%\begin{definition}
%    \label{def:polytope-bound}
%    Given a graph $G$ %and its induced algebraic load flow system $P_G$, as given 
%    %in \eqref{equ:powerflow-alg}, 
%    we define the \term{adjacency polytope bound} (or simply \term{AP bound}) 
%    for any $P_{G,Y}$ (with any choice of admittance matrix $Y$), denoted by 
%    $\mu_G$ to be the upper bound provided by Theorem \ref{thm:polytope-bound}. 
%    That is, $\mu_G \; := \; \nvol_{2n} (\nabla_G)$.
%\end{definition}

\begin{wrapfigure}[11]{r}{0.3\textwidth}\vspace{-0.35in}
    \centering
    \includegraphics[width=0.07\textwidth, height=0.2\textwidth]{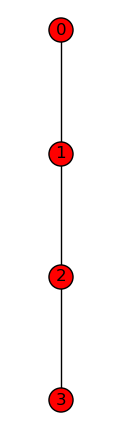}
    \caption{A path.}
    \label{fig:path-4}%\vspace{-0.35in}
\end{wrapfigure}

\begin{example}
    Consider, for example, the simple ``path graph'' $G = (\{0,1,2,3\},E)$ of 
    4 nodes as shown in Figure \ref{fig:path-4} where each node $i$ is connected 
    the next node $i+1$.
    Recall that we also require each node to have a loop to itself
    (to reflect the nonzero diagonal entries of the nodal admittance matrix $Y$),
    so the edges in the graph are
    \[
        E = \{
            (0,0), (0,1), (1,0), (1,1), (1,2), (2,1), (2,2), (2,3), (3,2), (3,3)
        \}.
    \]
    By Definition \ref{def:polytope}, the points in $\Gamma_G$ are therefore
    \begin{gather*}
        (\bolde_0,\bolde_0),
        (\bolde_0,\bolde_1),
        (\bolde_1,\bolde_0),
        (\bolde_1,\bolde_1),
        (\bolde_1,\bolde_2),
        (\bolde_2,\bolde_1),
        (\bolde_2,\bolde_2),
        (\bolde_2,\bolde_3),
        (\bolde_3,\bolde_2),
        (\bolde_3,\bolde_3).
    \end{gather*}
%    That is,
%    \begin{small}
%    \begin{equation*}
%        \Gamma_G =
%        \left\{
%        \begin{bmatrix}
%            0 \\ 0 \\ 0 \\ 0 \\ 0 \\ 0
%        \end{bmatrix},
%        \begin{bmatrix}
%            0 \\ 0 \\ 0 \\ 1 \\ 0 \\ 0
%        \end{bmatrix},
%        \begin{bmatrix}
%            1 \\ 0 \\ 0 \\ 0 \\ 0 \\ 0
%        \end{bmatrix},
%        \begin{bmatrix}
%            1 \\ 0 \\ 0 \\ 1 \\ 0 \\ 0
%        \end{bmatrix},
%        \begin{bmatrix}
%            1 \\ 0 \\ 0 \\ 0 \\ 1 \\ 0
%        \end{bmatrix},
%        \begin{bmatrix}
%            0 \\ 1 \\ 0 \\ 1 \\ 0 \\ 0
%        \end{bmatrix},
%        \begin{bmatrix}
%            0 \\ 1 \\ 0 \\ 0 \\ 1 \\ 0
%        \end{bmatrix},
%        \begin{bmatrix}
%            0 \\ 1 \\ 0 \\ 0 \\ 0 \\ 1
%        \end{bmatrix},
%        \begin{bmatrix}
%            0 \\ 0 \\ 1 \\ 0 \\ 1 \\ 0
%        \end{bmatrix},
%        \begin{bmatrix}
%            0 \\ 0 \\ 1 \\ 0 \\ 0 \\ 1
%        \end{bmatrix}
%        \right\}.
%    \end{equation*}
%    \end{small}
    Since $\boldzero$ is already included in $\Gamma_G$, the symmetric
    adjacency polytope $\nabla_G$ is precisely the convex hull of $\Gamma_G$.
    With programs for computing volume of convex polytopes to be listed in
    \S \ref{sec:computing-bound}, one can easily obtain that
    $\mu_G = \nvol_6 (\nabla_G) = \nvol_6 (\conv \Gamma_G) = 8$.
    That is, the algebraic load flow equation for such a path graph has at most
    8 isolated non-deficient complex solutions.
%    \begin{figure} %{r}{0.4\textwidth}
%        \centering
%        \includegraphics[width=0.5\textwidth]{gfx/path-4.png}
%        \caption{A graph with 4 nodes.}
%        \label{fig:path-4}
%    \end{figure}
\end{example}
\smallskip

%When compared to the BKK bound described in \S \ref{sec:bernshtein},
%the adjacency polytope bound is much easier to analyze.
%In particular, 
%Since the adjacency polytope bound is formulated in terms of volume,
It is quite easy to understand the situation under which the introduction of a 
new connection in the power network will change the AP bound:
%In this section we investigate how the AP bound (Definition
%\ref{def:polytope-bound}) changes as a connection or a new bus is introduced 
%into a power network.
%We start with some simple observations.
Since the AP bound is formulated in terms of the volume of
a polytope which is nondecreasing (i.e., it will either increase or remain
unchanged when new points are added), this upper bound must also be
nondecreasing with respect to the additions of new connections.

\medskip
\begin{theorem}
    \label{thm:nondecreasing}
    For a graph $G=(B,E)$ and two of its nodes $i$ and $j$ that are not 
    directly connected (i.e., $(i,j) \notin E$),
    let $G' = (B,E \cup \{(i,j)\})$ be the new graph constructed by adding 
    the edge between $i$ and $j$ to $G$. Then
    \begin{equation*}
        \mu_{G} \le \mu_{G'}.
    \end{equation*}
    Moreover, if $\{(\bolde_i,\bolde_j),(\bolde_j,\bolde_i)\} \subset \nabla_G$
    then
    \begin{equation*}
        \mu_{G} = \mu_{G'}.
    \end{equation*}
\end{theorem}
\medskip

\begin{proof}
    Recall that each edge in a graph contributes certain points (which may or 
    may not be vertices) in the construction of the symmetric adjacency polytope.
    Since the edges of $G$ is a subset of the edges of $G'$, we can see that
    $\nabla_G \subseteq \nabla_{G'}$ with the equality hold precisely when the
    points contributed by $(i,j)$ are already inside $\nabla_G$.
    With these observations in mind, both parts of this theorem are direct 
    implications of normalized volume being nondecreasing.
    %    This property itself is not immediate obvious from the definition of
    %    mixed volume. Here we refer to the approach of mixed cells 
    %    \cite{huber_polyhedral_1995,chen_theoretical_2014} for a geometric explanation.
\end{proof}
\medskip

Based on this observation, it can be shown that the AP bound is never more than 
the BBLSY bound.
This is essentially our alternative proof of the BBLSY bound:

\medskip
\begin{theorem}
    For a graph $G=(B,E)$, let $|B| = n + 1$. Then, $\mu_G \le \binom{2n}{n}$
\end{theorem}
\medskip

\begin{proof}
    Fixing the set of buses, Theorem \ref{thm:nondecreasing} states that the 
    AP bound is nondecreasing as new edges are added to the graph.
    Consequently, the AP bound for any network constructed from this set of buses
    is bounded above by the AP bound for the graph with most edges,
    that is, a complete graph.
    It is easy to verify that for a complete graph $G = (B,E)$ 
    (with loops), $\nabla_G \subseteq (\conv A) + (\conv B)$ where
    \begin{align}
    A &= \left\{ \bolde_0, \bolde_{  1}, \bolde_{  2},\ldots\,\bolde_{ n} \right\} &
    B &= \left\{ \bolde_0, \bolde_{n+1}, \bolde_{n+2},\ldots\,\bolde_{2n} \right\}
    \end{align}
    and $(\conv A)+(\conv B)$ denotes the Minkowski sum of the two polytopes
    $(\conv A)$ and $(\conv B)$.
    Note that both $(\conv A)$ and $(\conv B)$ are $n$-dimensional.
    Then by multi-linearity of mixed volume with respect to Minkowski sum,
    \begin{equation*}
    \begin{aligned}
    \nvol((\conv A) + (\conv B)) &=
    \sum_{k=0}^{2n} \binom{2n}{k} MV((conv A)^{(k)},(conv B)^{(2n-k)}) \\ &=
    \binom{2n}{n}  MV((conv A)^{(n)},(conv B)^{(n)})\\ &=
    \binom{2n}{n}
    \end{aligned}
    \end{equation*}
\end{proof}
\medskip

We conclude this section with a reiteration of the various root counts
involved in the discussion.
Recall that for a power network $G = (B,E)$, the number of isolated 
non-deficient solutions of the original (non-algebraic) load flow equation
\eqref{equ:Sj} (physical solutions), the number of isolated complex solutions of 
the algebraic load flow equation \eqref{equ:powerflow-alg} in $\Cstar^{2n}$, 
the BKK bound provided by Theorem \ref{thm:bernshtein}, 
the AP bound given in Theorem \ref{thm:polytope-bound}, 
the BBLSY bound given in \cite{baillieul1982,li_numerical_1987},
and the CB number are related as follows:
%\begin{small}
\begin{equation*}
    \parbox{9ex}{\centering Physical \\ solution \\ count} \;\le\;
    \parbox{9ex}{\centering Complex \\ solution \\ count} \;\le\;
    \parbox{7ex}{\centering BKK \\ bound} \;\le\;
    \parbox{7ex}{\centering AP \\ bound} \;\le\;
    \parbox{8ex}{\centering BBLSY \\ bound} \;\le\;
    \parbox{9ex}{\centering CB \\ bound}
\end{equation*}
\subsection{Computing BKK and AP bounds} \label{sec:computing-bound}

The BKK bound, can be computed by using
efficient software packages such as 
\tech{DEMiCs} \cite{mizutani_demics_2008},
\tech{Gfan} \cite{jensen_presentation_2006},
\tech{MixedVol} \cite{gao_algorithm_2005}, 
\tech{MixedVol-2.0} \cite{lee_mixed_2011}.
For larger power networks involving many buses, the induced 
algebraic load flow equation may contain a large number of terms.
Parallel computing technology will be essential in the computation of the 
BKK bounds for such large scale polynomial systems.
\tech{MixedVol-3} \cite{chen_mixed_2014,chen_mixed}
(with an improved version integrated in \tech{Hom4PS-3} \cite{chen_hom4ps3_2014})
is capable of computing the BKK bound for large polynomial systems in parallel
on a wide range of hardware architectures including multi-core systems,
NUMA systems, and computer clusters.
As noted in Remark \ref{rmk:union-benefit}, however, there is built-in level of
redundancy in the Newton polytopes (see Appendix \ref{sec:mvol}) of the 
algebraic load flow equations. 
The formulation of the AP bound takes advantage of this natural 
redundancy and can generally be computed much more easily than the BKK bound for
larger power networks.
Since the AP bound is formulated in terms of the volume of a
convex polytope (the symmetric AP), any software that can
compute such volume exactly can be used to provide this bound.
A survey on the various algorithms for exact volume computation can be found in
\cite{bueler_exact_2000}.
A particularly versatile software package that utilizes a wide range of algorithms
is the \tech{Vinci} program \cite{enge_vinci}.

%===============================================================================
\subsection{Homotopy methods for solving load flow equations}
\label{sec:sol_count_homotopy}
The previous sections described two upper bounds --- the BKK bound and the AP bound
--- for the number of isolated non-deficient complex solutions of the 
algebraic load flow equation induced by a given network topology.
It is worth reiterating that the BKK bound is more than just an \emph{upper bound}:
As shown in Theorem \ref{thm:attainable}, it is actually the \emph{maximum} 
for the given network topology in the sense that there always exits \emph{some} 
choice of the admittance matrix $Y$ (and $S$) for which the total number of
isolated non-deficient complex solutions is exactly the BKK bound.
While the AP bound, in general, may be larger than the BKK bound, as we shall
show in \S \ref{sec:cases}, the two coincide for all the networks we have 
investigated in the present work.
The family of numerical methods known as 
homotopy methods have been proved to be a robust and efficient approach for 
solving algebraic load flow equations.
One great strength of these methods lies in the pleasantly parallel nature:
in principle, each solution can be computed independently.
This feature is of particular importance in dealing with larger power networks 
containing many buses and connections (hence more complicated equations).
It is therefore a natural question to ask: is there a homotopy method that
can solve the algebraic load flow equation by tracking BKK bound number of
homotopy paths?
This section establishes the answer in the affirmative.

This homotopy method is the \emph{polyhedral homotopy} method developed in
\cite{huber_polyhedral_1995}. Here we briefly state the construction.
We choose a pair of random (rational) numbers $\omega_{ij}$ and $\omega_{ij}'$
for each $i,j=1,\dots,n$.
%For each term $v_i u_j$ appeared in \eqref{equ:powerflow-alg}, in the first $n$
%equations we choose a random (rational) number he lifting for the term   will be denoted by $\omega_{ij}$ while that of $v_i u_j$ in the
%last $n$ equations will be denoted by $\omega_{ij}'$.
With these we define the \emph{homotopy function}
\begin{equation}
    H_{G,Y} (v_1,\dots,v_n,u_1,\dots,u_n,t) = 
    \begin{cases}
        \sum_{k=0}^{n} Y_{1k}\, (t) v_1 u_k t^{\omega_{1k}} - S_1(t)   \\
        \quad \vdots \\
        \sum_{k=0}^{n} Y_{nk}   (t) v_n u_k t^{\omega_{nk}} - S_n(t)   \\[1.5ex]
        \sum_{k=0}^{n} Y_{1k}'\,(t) u_1 v_k t^{\omega_{k1}'} - S_1'(t) \\
        \quad \vdots \\
        \sum_{k=0}^{n} Y_{nk}'  (t) u_n v_k t^{\omega_{kn}'} - S_n'(t) .
    \end{cases}
    \label{equ:polyhedral}
\end{equation}
where
\begin{align*}
    Y_{ij} (t) &= (1-t) Z_{ij}  + t Y_{ij}^* & S_{i} (t) &= (1-t)W_i  + S_i \\
    Y_{ij}'(t) &= (1-t) Z_{ij}' + t Y_{ij}   & S_{i}'(t) &= (1-t)W_i' + S_i^*
\end{align*}
and $(Z_{ij})$ and $(Z_{ij}')$ are randomly chosen complex matrices of the same 
size as $Y$ and $W = (W_i)$ and $W' = (W_i')$ are two random complex vectors 
in $\C^n$.

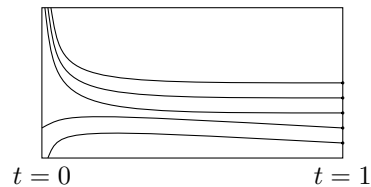
\begin{wrapfigure}{r}{0.4\columnwidth}
    \centering
    \begin{tikzpicture}[scale=0.4]
    \draw (0,0) rectangle (10,5);
    \draw (0,0) node[below] {$t=0$};
    \draw (10,0) node[below] {$t=1$};
    \draw [] (10,0.5) circle(1pt) .. controls (0.6,1.0) .. (0.2,0);
    \draw [] (10,1.0) circle(1pt) .. controls (0.9,1.5) .. (0,1);
    \draw [] (10,1.5) circle(1pt) .. controls (0.5,1.5) .. (0.1,5);
    \draw [] (10,2.0) circle(1pt) .. controls (0.6,2.0) .. (0.2,5);
    \draw [] (10,2.5) circle(1pt) .. controls (0.7,2.5) .. (0.3,5);
    %\draw [thick,-latex] (0,3.0) .. controls (9,3.0) .. (9.6,5);
    %\draw [thick,->] (0,3.5) .. controls (9,3.5) .. (9.5,5);
    %\draw [thick,->] (0,4.0) .. controls (9,4.0) .. (9.4,5);
    \end{tikzpicture}
    \caption{The solution paths defined by $H_{G,Y}(\boldv,\boldu,t) = \boldzero$.}
    \label{fig:polyhedral-start}\vspace{-0.3in}
\end{wrapfigure}

Clearly $H_{G,Y}(\boldv,\boldu,1) \equiv P_{G,Y}(\boldv,\boldu)$.
For generic choice of $Z$, $Z'$, $W$, $W'$, $\omega$ and $\omega'$,
it can be shown that fixing any $t \in (0,1)$, the non-deficient solutions of
$H_{G,Y}(\boldv,\boldu,t) = \boldzero$ are all isolated and the total number is
exactly the BKK bound.
Moreover, as $t$ varies within $(0,1)$, the corresponding solutions of 
$H_{G,Y}(\boldv,\boldu,t) = \boldzero$ also vary smoothly forming smooth
\emph{solution paths} that collectively reach all the desired solutions of
the algebraic load flow equation $P_{G,Y}(\boldv,\boldu) = \boldzero$.
See Figure \ref{fig:polyhedral-start} for a schematic illustration.
Thus, once the ``starting points'' of each solution path at $t=0$ are found,
standard numerical continuation techniques can be used to track the
solution paths and reach \emph{all} the isolated non-deficient complex solutions.

An apparent difficulty is the identifications of the ``starting points''.
After all, at $t=0$, $H_{G,Y}(\boldu,\boldv,t)$ either becomes constant
or undefined since each of its nonconstant term contains a power of $t$.
This, fortunately, is surmounted via a device known as \emph{mixed cells}
which are themselves the by-product of the computation of the BKK bound.
The technical detail is outside the scope of the present article,
we therefore refer to standard references such as 
\cite{huber_polyhedral_1995,li_numerical_2003}.
The application of the polyhedral homotopy to load flow equations will be
explored in future works, here we simply emphasize that with the
polyhedral homotopy method, the number of solution paths one needs to track
is precisely the BKK bound \textit{of the particular polynomial formulation of the load flow equations
that we have used (eq. (\ref{equ:powerflow-alg}))} --- the best possible network topology-based
complex solution bound.

\FloatBarrier
%===============================================================================
\section{Solution bound for certain power networks} \label{sec:cases}

In this section we provide concrete computation results for the BKK bound and 
the AP bound described above when applied to specific graphs.
We should reiterate that we require all the graphs to have self-loops for each 
node to reflect the nonzero diagonal entries of the nodal admittance matrix $Y$.
In all cases, the BKK and the AP bounds are computed via \tech{MixedVol-3} 
\cite{chen_mixed_2014,chen_mixed}.
Complex solutions count of specific load flow systems are computed by actually
solving the systems via \tech{Hom4PS-3} \cite{chen_hom4ps3,chen_hom4ps3_2014}
which implements the \emph{polyhedral homotopy method} described above.

For comparison, in each case we show the ``5-way comparison'' among the 
complex solution count, BKK bound, AP bound, the BBLSY Bound, and the 
CB number.
Here, ``complex solution count'' refers to the number of isolated complex 
solutions of \eqref{equ:powerflow-alg} in $\Cstar^{2n}$ where $n$ is number 
of non-reference nodes.
Since the actual complex solution count may depend on the coefficients
($Y_{ij}$ and $S_i$ in \eqref{equ:powerflow-alg}), they are computed based on 
a sample of randomly chosen set of coefficients for each graph.

\subsection{Paths}

We first consider an extremely sparsely connected family of graphs -- path graphs.
In a path graph $G = (B,E)$ with nodes $0,1,\dots,n$ where $n = |B|-1$
(the number of non-reference nodes), a node $i$ is connected to two of its neighbors 
$i+1$ and $i-1$, with the exceptions that node $0$ is only connected to node 1
and the last node, node $n$, is only connected to node $n-1$.

Table \ref{tab:path} shows the ``5-way comparison'' described above.
It is important to note that in all cases computed (100 in total, with 10 random 
admittance matrices chosen for each $|B|$), the actual complex solution count,
the BKK bound, and the AP bound are exactly the same.

%In \cite{baillieul1982}, based on Morse inequalities, it was shown that a 
%lower bound on the number of load flow solutions was $2^{n-1}$. 
%We have shown that the system on the path graph certainly attains this lower bound.
%It would be interesting to come up with other network topologies that also attains this lower bound.

\begin{table}[h!]
    \centering
    \scriptsize
    \begin{tabular}{lrrrrrrrrrrr}
    	\toprule
    	$|B|$     & 2 &  3 &  4 &   5 &    6 &    7 &     8 &     9 &     10 &      11 &      12 \\ \midrule
    	Solutions & 2 &  4 &  8 &  16 &   32 &   64 &   128 &   256 &    512 &    1024 &    2048 \\
    	BKK       & 2 &  4 &  8 &  16 &   32 &   64 &   128 &   256 &    512 &    1024 &    2048 \\
    	AP        & 2 &  4 &  8 &  16 &   32 &   64 &   128 &   256 &    512 &    1024 &    2048 \\
    	BBLSY        & 2 &  6 & 20 &  70 &  252 &  924 &  3432 & 12870 &  48620 &  184756 &  705432 \\
    	CB        & 4 & 16 & 64 & 256 & 1024 & 4096 & 16384 & 65536 & 262144 & 1048576 & 4194304 \\ \bottomrule
    \end{tabular}
    \caption{
        The ``5-way comparison'' of the solution bounds for path graphs of size 2--12.
    }
    \label{tab:path}\vspace{-0.2in}
\end{table}

%\textit{@Tianran: Any further concrete way to 'prove' this $2^{n-1}$ number? DanM specifically suggested this.}

\FloatBarrier
%===============================================================================
%\vspace{-0.3in}
\subsection{Rings}
\begin{figure}[h!]
    \centering
    \begin{subfigure}[t]{0.4\textwidth}
 \centering
    \includegraphics[width=0.7\columnwidth]{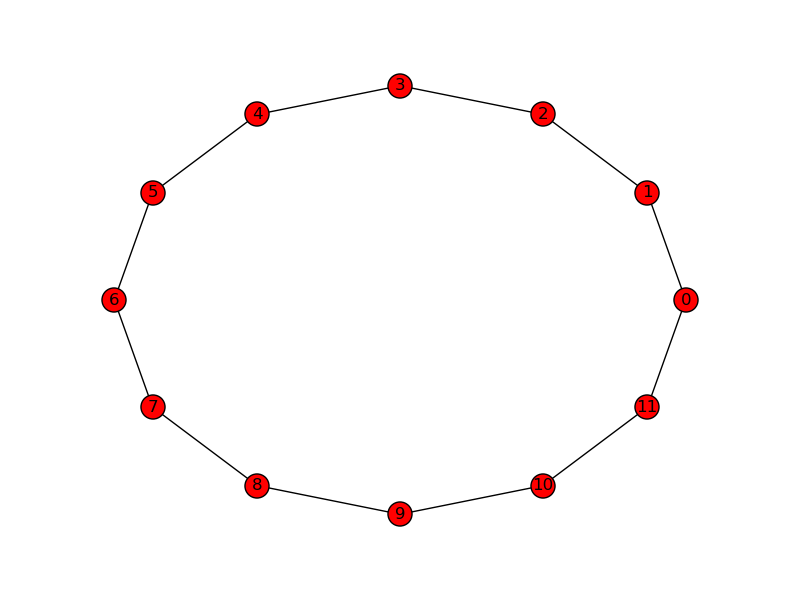}
    \caption{A ring.}
    \label{fig:ring}%\vspace{-0.3in}
    \end{subfigure}%
    \begin{subfigure}[t]{0.4\textwidth}
        \centering
\includegraphics[width=0.7\columnwidth]{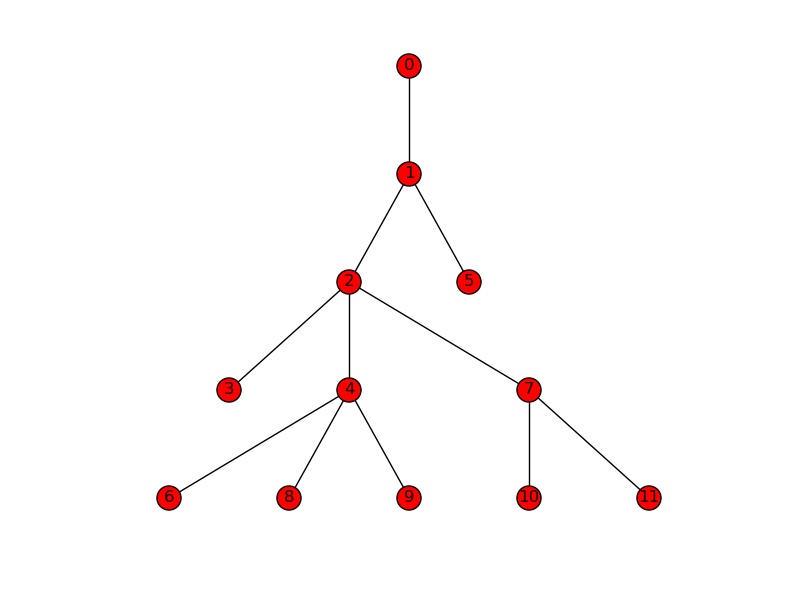}
    \caption{A random tree}
    \label{fig:rand-tree}
    \end{subfigure}%
\end{figure}

By joining the two ends of a path graph studied above, 
a ``ring graph'' is obtained.
More precisely, in a ring $G = (B,E)$ with nodes $0,1,\dots,n$ where $n = |B|-1$, 
a node $i$ is connected to two of its neighbors $i+1$ and $i-1$, and node reference
node, node $0$ is connected to node $n$. 
See, for example, the graph shown in Figure \ref{fig:ring}.
Table \ref{tab:ring} shows the similar ``5-way comparison''.
%to the table above: the 5-way
%comparison among the minimum number of isolated non-deficient complex solutions 
%in 10 random cases, the BKK bound, the AP bound, the BB bound, 
%and the CBB.
Again, in all cases computed (100 in total, with 10 random admittance matrices 
chosen for each $|B|$), the actual complex solution count, the BKK bound, 
and the AP bound are all exactly the same.

\begin{comment}
%\begin{wrapfigure}{l}{0.35\textwidth}
 \begin{figure}[h!]
 \centering
    \includegraphics[width=0.3\columnwidth]{gfx/ring-12.png}
    \caption{A ring.}
    \label{fig:ring}%\vspace{-0.3in}
\end{figure}
    %\end{wrapfigure}%
\end{comment}

%From the data in Table \ref{tab:ring}, the BKK and AP (and the number of complex solutions, at least for the 
%generic complex parameter values) appear to be $n 2^{n-2}$. \textit{@Tianran: can we prove this at all? If not,
%we can probably mention it as a conjecture for someone (may be the future us?!) to prove it?}

%\begin{figure}[h]
%    \centering
%    \includegraphics[width=0.6\textwidth]{gfx/ring-12.png}
%    \caption{A ring consisting of 12 nodes.}
%    \label{fig:ring}
%\end{figure}

\begin{table}[h!]
    \centering
    \scriptsize
    \begin{tabular}{lrrrrrrrrrrr}
    	\toprule
    	$|B|$     & 2 &  3 &  4 &   5 &    6 &    7 &     8 &     9 &     10 &      11 &      12 \\ \midrule
    	Solutions & 2 &  6 & 16 &  40 &   96 &  224 &   512 &  1152 &   2560 &    5632 &   12288 \\
    	BKK       & 2 &  6 & 16 &  40 &   96 &  224 &   512 &  1152 &   2560 &    5632 &   12288 \\
    	AP        & 2 &  6 & 16 &  40 &   96 &  224 &   512 &  1152 &   2560 &    5632 &   12288 \\
    	BBLSY        & 2 &  6 & 20 &  70 &  252 &  924 &  3432 & 12870 &  48620 &  184756 &  705432 \\
    	CB        & 4 & 16 & 64 & 256 & 1024 & 4096 & 16384 & 65536 & 262144 & 1048576 & 4194304 \\ \bottomrule
    \end{tabular}
    \caption{
        The ``5-way comparison'' of solution bounds for ring graphs of size 2--12.
        }
    \label{tab:ring}
\end{table}

\subsection{Random trees}
\begin{comment}
\begin{wrapfigure}[11]{r}{0.35\columnwidth}\vspace{-0.32in}
    \includegraphics[width=0.36\columnwidth]{gfx/rand-tree-12.png}
    \caption{A random tree}
    \label{fig:rand-tree}
\end{wrapfigure}
\end{comment}
Since distribution networks often have structures of trees, randomly generated 
trees may actually resemble some realistic power networks.
To get a reasonable sample size, tree graphs with $|B|=3,\dots,13$ are studied. 
For each of the sizes 10 trees are randomly generated (see, for example, the
trees shown in Figure \ref{fig:rand-tree}).
For each of these trees, 100 set of random $Y$ matrices are assigned.
Therefore for each given size $|B|$, 1000 polynomial systems are generated
(a total of 10000 systems for $|B|=3,\dots,13$).
%They are solved by \tech{Hom4PS-3} \cite{chen_hom4ps3,chen_hom4ps3_2014}, and 
Table \ref{tab:rand-tree} shows the similar ``5-way comparison'' used above.
%Remarkably, the actual number of complex solutions of these 10000 systems in 
%$\Cstar^{2(|B|-1)}$ all agree with the BKK and AP bounds as shown in. 
%The BKK and AP bounds and the 
%number of complex solutions at generic complex parameter points also are clearly $2^{n-1}$. 
%It must also be emphasized that the BKK (and AP) bounds is dramatically smaller 
%compared to the BB or CB bound. 
%For example, for $n=13$, the BKK bound is ~$1/660$th of the BB bound.

\begin{table}[h!]
    \centering
    \tiny
    \begin{tabular}{lrrrrrrrrrrr}
    	\toprule
    	$|B|$     &  3 &  4 &   5 &    6 &    7 &     8 &     9 &     10 &      11 &      12 &       13 \\ \midrule
    	Solutions &  4 &  8 &  16 &   32 &   64 &   128 &   256 &    512 &    1024 &    2048 &     4096 \\
    	BKK       &  4 &  8 &  16 &   32 &   64 &   128 &   256 &    512 &    1024 &    2048 &     4096 \\
    	AP        &  4 &  8 &  16 &   32 &   64 &   128 &   256 &    512 &    1024 &    2048 &     4096 \\
    	BBLSY     &  6 & 20 &  70 &  252 &  924 &  3432 & 12870 &  48620 &  184756 &  705432 &  2704156 \\
    	CB        & 16 & 64 & 256 & 1024 & 4096 & 16384 & 65536 & 262144 & 1048576 & 4194304 & 16777216 \\ \bottomrule
    \end{tabular}
    \caption{
        The ``5-way comparison'' for randomly generated trees of given sizes.
    }
    \label{tab:rand-tree}
\end{table}

%\begin{figure}[h]
%    \centering
%    \begin{subfigure}[t]{0.48\columnwidth}
%        \includegraphics[width=\columnwidth]{gfx/rand-tree-8.png}
%        %\caption{$|B|=8$}
%    \end{subfigure}%
%    \begin{subfigure}[t]{0.48\columnwidth}
%        \includegraphics[width=\columnwidth]{gfx/rand-tree-12.png}
%        %\caption{$|B|=12$}
%    \end{subfigure}
%    \caption{Two randomly generated trees.}
%    \label{fig:rand-tree}
%\end{figure}

%\FloatBarrier
%===============================================================================
\subsection{Clusters}

Real power networks generally exhibit certain level of ``clustering'', that is, 
certain subset of buses are densely connected while on a larger scale,
the connections among such subsets are sparse.
%For example, within a city, the buses can be densely connected with one another in 
%the sense that most buses are connected to most other buses while between cities
%the connections are sparse as only very few long distance power lines make the 
%connection.
%In Section \ref{sec:dense}, we have discussed the well known cases of completely 
%connected networks. 
Here for simplicity, we focus on the most extreme cases where a 
larger network is created by joining completely connected subnetworks.
This section lists concrete results of the BKK and AP bounds for certain 
simple types of clusters.
For comparison, in each case we only show the ``3-way comparison'' among the 
BKK bound, AP bound, and the BBLSY bound (due to the large amount of data).
%Here, ``complex solution count'' refers to the number of isolated complex 
%solutions of \eqref{equ:powerflow-alg} in $\Cstar^{2n}$ where $n$ is number 
%of non-reference nodes.
%Since the actual complex solution count may depend on the coefficients
%($Y_{ij}$ and $S_i$ in \eqref{equ:powerflow-alg}), they are computed based on 
%a sample of randomly chosen set of coefficients for each graph.

%\FloatBarrier
%===============================================================================
\subsubsection{Subnetworks sharing nodes}

\begin{figure}[h!]
    \centering
    \begin{subfigure}[t]{0.49\textwidth}
        \centering
        \includegraphics[width=0.8\columnwidth]{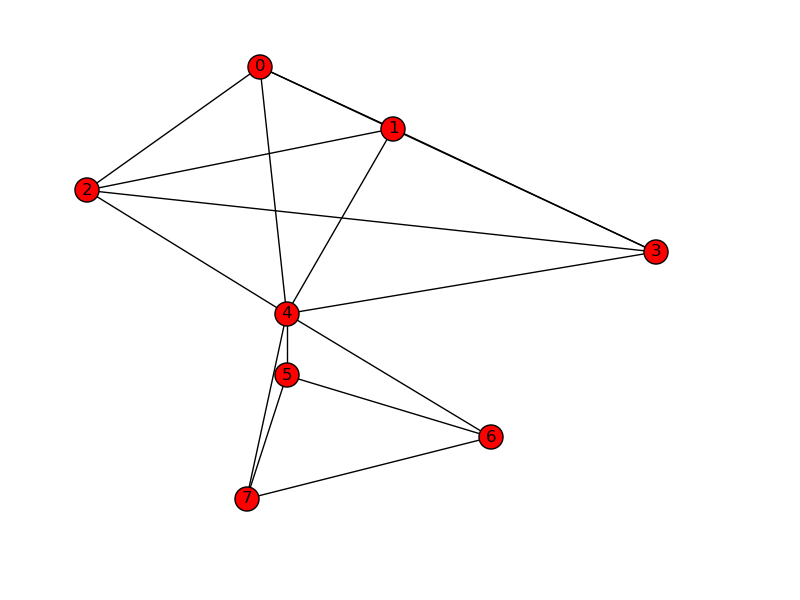}
        \caption{Sharing one node}
        \label{fig:2cluster-share1}
    \end{subfigure}%
    \begin{subfigure}[t]{0.49\textwidth}
        \centering
        \includegraphics[width=0.8\columnwidth]{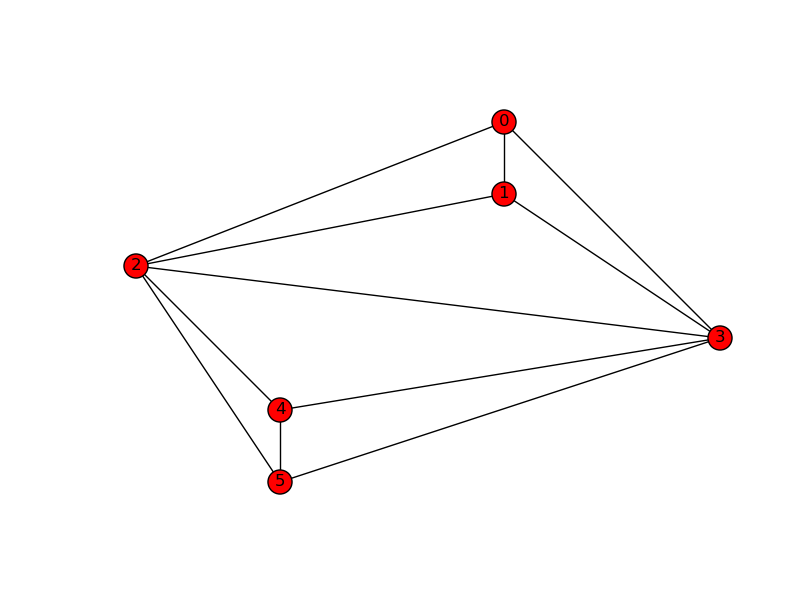}
        \caption{Sharing two nodes}
        \label{fig:2cluster-share2}
    \end{subfigure}%
%    \begin{subfigure}[t]{0.33\textwidth}
%        \includegraphics[width=0.9\columnwidth]{gfx/2cluster-share3-5-6.png}
%        \caption{Sharing three nodes}
%        \label{fig:2cluster-share3}
%    \end{subfigure}%
    \caption{Completely connected subgraphs sharing common nodes.}
    \label{fig:2cluster-share}
\end{figure}

%in the subset is connected directly to every other node in the subset.
%Such subsets are known as \emph{cliques} in graph theory.
%In particular, we shall investigate networks formed by connecting several
%cliques by certain edges or common nodes.
%In graph theory, a (maximal) clique is a subgraph in which every node is
%connected to every other node and the inclusion of any other node will
%violate this property.
%In other words, in a network, a clique is a (most) densely connected subnetwork
%that cannot be extended any further.
%A clique is a convenient unit in studying the large scale structures in a cluster.
%We start with the cases of two cliques sharing a single common non-reference bus.

First, these subnetworks can be jointed by sharing common nodes.
See, for example, the networks shown in Figure \ref{fig:2cluster-share}.
Table \ref{tab:2cluster-share1-all} shows the ``3-way comparison'' for cases
where two completely connected subnetworks share 
\textit{a single (non-reference) common bus}.
These cases have been studied in \cite{guo1990}.
Our computational results agree with their assertion.
%Note that the complex solution count, the BKK bound, and the AP bound are 
%exactly the same for the cases shown in the table.
Table \ref{tab:2cluster-share2-all} shows the similar ``3-way comparison'' 
for cases where two completely connected subnetworks share 
\textit{two (non-reference) common bus}.
These cases have been extensively studied in \cite{molzahn-mehta-matt2015}
via numerical methods.
The results and conjectures in that work are precisely reproduced
by our computation.
%The results for for cases where two completely connected subnetworks share 
%\textit{three (non-reference) common bus} is shown in 
%Table \ref{tab:2cluster-share3-all}.
For larger networks, the AP bounds are generally much easier to compute
than the BKK bound using existing implementations. 
In Table \ref{tab:2cluster-share1} and \ref{tab:2cluster-share2}, 
we show the AP bounds for the above referenced clusters.

\begin{table}[h!]
    \centering
    \tiny
    \begin{tabular}{rrrrrr}
    	\toprule
    	$c_1 \backslash c_2$ &           2 &              3 &               4 &                 5 &                  6 \\ \midrule
    	                   2 &       4/4/6 &       12/12/20 &        40/40/70 &       140/140/252 &        504/504/924 \\
    	                   3 &    12/12/20 &       36/36/70 &     120/120/252 &       420/420/924 &     1512/1512/3432 \\
    	                   4 &    40/40/70 &    120/120/252 &     400/400/924 &    1400/1400/3432 &    5040/5040/12870 \\
    	                   5 & 140/140/252 &    420/420/924 &  1400/1400/3432 &   4900/4900/12870 &  17640/17640/48620 \\
    	                   6 & 504/504/924 & 1512/1512/3432 & 5040/5040/12870 & 17640/17640/48620 & 63504/63504/184756 \\ \bottomrule
    \end{tabular}
    \caption{
        The ``3-way comparison'' %among the BKK, AP, and the BB bound 
        for two completely connected subnetworks of sizes $c_1$ and $c_2$ 
        respectively sharing a \textbf{single} common non-reference nodes.
    }
    \label{tab:2cluster-share1-all}
\end{table}

\begin{table}[h!]
    \centering
    \tiny
    \begin{tabular}{rrrrrr}
    	\toprule
    	$c_1 \backslash c_2$ &           2 &           3 &              4 &               5 &                 6 \\ \midrule
    	                   2 &       2/2/2 &       6/6/6 &       20/20/20 &        70/70/70 &       252/252/252 \\
    	                   3 &       6/6/6 &    18/18/20 &       60/60/70 &     210/210/252 &       756/756/924 \\
    	                   4 &    20/20/20 &    60/60/70 &    200/200/252 &     700/700/924 &    2520/2520/3432 \\
    	                   5 &    70/70/70 & 210/210/252 &    700/700/924 &  2450/2450/3432 &   8820/8820/12870 \\
    	                   6 & 252/252/252 & 756/756/924 & 2520/2520/3432 & 8820/8820/12870 & 31752/31752/48620 \\ \bottomrule
    \end{tabular}
    \caption{
        The ``3-way comparison'' for two completely connected subnetworks of 
        sizes $c_1$ and $c_2$ respectively sharing 
        \textbf{two} common non-reference nodes.
    }
    \label{tab:2cluster-share2-all}
\end{table}

%\begin{table}[h]
%    \centering
%    \tiny
%    \begin{tabular}{rrrrr}
%    	\toprule
%    	$c_1 \backslash c_2$ &           3 &           4 &              5 &                 6 \\ \midrule
%    	                   3 &       6/6/6 &    20/20/20 &       70/70/70 &       252/252/252 \\
%    	                   4 &    20/20/20 &    68/68/70 &    240/240/252 &       868/868/924 \\
%    	                   5 &    70/70/70 & 240/240/252 &    850/850/924 &    3080/3080/3432 \\
%    	                   6 & 252/252/252 & 868/868/924 & 3080/3080/3432 & 11172/11172/12870 \\ \bottomrule
%    \end{tabular}
%    \caption{
%        The ``3-way comparison'' for two completely connected subnetworks 
%        of sizes $c_1$ and $c_2$ respectively sharing 
%        \textbf{three} common non-reference nodes.
%    }
%    \label{tab:2cluster-share3-all}
%\end{table}

\begin{table}[h!]
    \centering
    \tiny
    \begin{tabular}{rrrrrrrr}
    	\toprule
    	$c_1 \backslash c_2$ &    2 &     3 &     4 &      5 &      6 &       7 &        8 \\ \midrule
    	                   2 &    4 &    12 &    40 &    140 &    504 &    1848 &     6864 \\
    	                   3 &   12 &    36 &   120 &    420 &   1512 &    5544 &    20592 \\
    	                   4 &   40 &   120 &   400 &   1400 &   5040 &   18480 &    68640 \\
    	                   5 &  140 &   420 &  1400 &   4900 &  17640 &   64680 &   240240 \\
    	                   6 &  504 &  1512 &  5040 &  17640 &  63504 &  232848 &   864864 \\
    	                   7 & 1848 &  5544 & 18480 &  64680 & 232848 &  853776 &  3171168 \\
    	                   8 & 6864 & 20592 & 68640 & 240240 & 864864 & 3171168 & 11778624 \\ \bottomrule
    \end{tabular}
    \caption{
        The AP bounds for graphs consisting of two cliques of size 
        $c_1$ and $c_2$ respectively sharing a single common non-reference node.
        }
    \label{tab:2cluster-share1}
\end{table}

\begin{table}[h!]
    \centering
    \tiny
    \begin{tabular}{rrrrrrrr}
        \toprule
        $c_1 \backslash c_2$ &    2 &     3 &     4 &      5 &      6 &       7 &       8 \\ \midrule
        2 &    2 &     6 &    20 &     70 &    252 &     924 &    3432 \\
        3 &    6 &    18 &    60 &    210 &    756 &    2772 &   10296 \\
        4 &   20 &    60 &   200 &    700 &   2520 &    9240 &   34320 \\
        5 &   70 &   210 &   700 &   2450 &   8820 &   32340 &  120120 \\
        6 &  252 &   756 &  2520 &   8820 &  31752 &  116424 &  432432 \\
        7 &  924 &  2772 &  9240 &  32340 & 116424 &  426888 & 1585584 \\
        8 & 3432 & 10296 & 34320 & 120120 & 432432 & 1585584 & 5889312 \\ \bottomrule
    \end{tabular}
    \caption{
        The AP bounds for graphs consisting of two cliques, of size $c_1$ and
        $c_2$ respectively sharing two common non-reference nodes.
    }
    \label{tab:2cluster-share2}
\end{table}

\subsubsection{Completely connected subnetworks connected by edges}

\begin{figure}[h]
    \centering
    \begin{subfigure}[t]{0.45\textwidth}
        \centering
        \includegraphics[width=0.95\columnwidth]{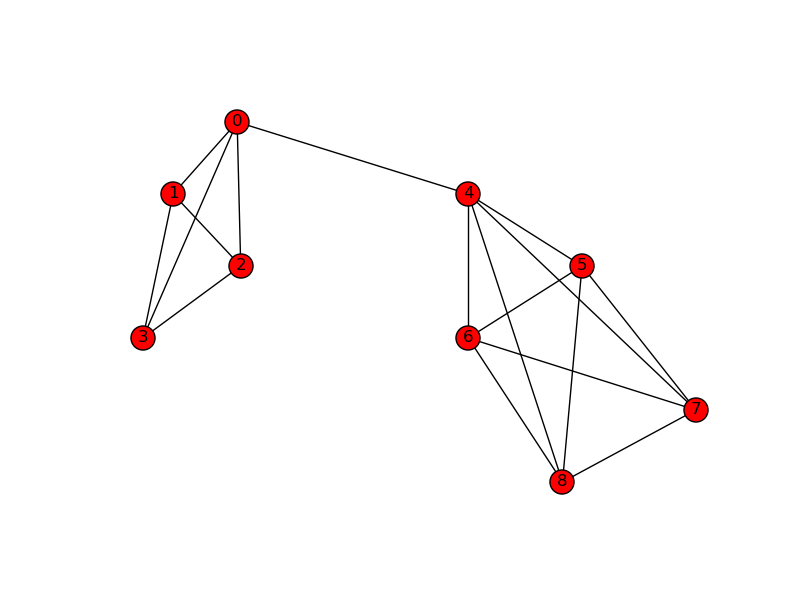}
        \caption{Two completely connected subnetworks jointed via a single edge.}
        \label{fig:2cluster-edge}
    \end{subfigure}
    \hfill
    \begin{subfigure}[t]{0.45\textwidth}
        \centering
        \includegraphics[width=0.99\columnwidth]{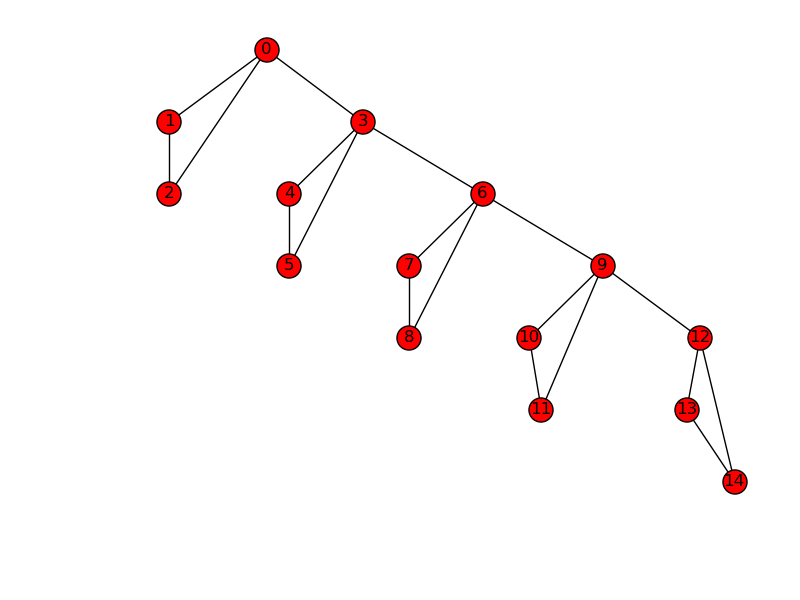}
        \caption{A network containing several subnetworks linked together via edges.}
        \label{fig:mcluster}
    \end{subfigure}
    \caption{Completely connected subnetworks connected by edges}
\end{figure}

We now consider cases where completely connected subnetworks are connected by edges.
For example, Figure \ref{fig:2cluster-edge} shows a network that consists of 
two cliques of size four and five respectively connected by a single edge.
Table \ref{tab:2cluster} shows the AP bounds for networks created from joining
two completely connected subnetworks by a \textit{single} edge.

%\begin{table}[h]
%    \centering
%    \tiny
%    \begin{tabular}{rrrrrr}
%    	\toprule
%    	$c_1 \backslash c_2$ &                   2 &                    3 &                       4 &                        5 &                           6 \\ \midrule
%    	                   2 &            8/8/8/20 &          24/24/24/70 &            80/80/80/252 &          280/280/280/924 &         1008/1008/1008/3432 \\
%    	                   3 &         24/24/24/70 &         72/72/72/252 &         240/240/240/924 &         840/840/840/3432 &        3024/3024/3024/12870 \\
%    	                   4 &        80/80/80/252 &      240/240/240/924 &        800/800/800/3432 &     2800/2800/2800/12870 &     10080/10080/10080/48620 \\
%    	                   5 &     280/280/280/924 &     840/840/840/3432 &    2800/2800/2800/12870 &     9800/9800/9800/48620 &    35280/35280/35280/184756 \\
%    	                   6 & 1008/1008/1008/3432 & 3024/3024/3024/12870 & 10080/10080/10080/48620 & 35280/35280/35280/184756 & 127008/127008/127008/705432 \\ \bottomrule
%    \end{tabular}
%    \caption{
%        The 4-way comparison among the complex solution count, BKK bound,
%        AP bound, and the BB bound for graphs consisting of 
%        two cliques connected by a single edge.
%        The actual complex solution count is computed based on 10 randomly 
%        chosen admittance matrix for each case.
%    }
%    \label{tab:2cluster-all}
%\end{table}

\begin{table}[h!]
    \centering
    \tiny
    \begin{tabular}{crrrrrrrrrr}
    	\toprule
    	$c_1\backslash c_2$ &     1 &      2 &      3 &       4 &       5 &       6 &       7 &       8 &       9 &      10 \\ \midrule
    	         1          &       &      4 &     12 &      40 &     140 &     504 &    1848 &    6864 &   25740 &   97240 \\
    	         2          &     4 &      8 &     24 &      80 &     280 &    1008 &    3696 &   13728 &   51480 &  194480 \\
    	         3          &    12 &     24 &     72 &     240 &     840 &    3024 &   11088 &   41184 &  154440 &  583440 \\
    	         4          &    40 &     80 &    240 &     800 &    2800 &   10080 &   36960 &  137280 &  514800 & 1944800 \\
    	         5          &   140 &    280 &    840 &    2800 &    9800 &   35280 &  129360 &  480480 & 1801800 &  \\
    	         6          &   504 &   1008 &   3024 &   10080 &   35280 &  127008 &  465696 & 1729728 &         &  \\
    	         7          &  1848 &   3696 &  11088 &   36960 &  129360 &  465696 & 1707552 &         &         &  \\
    	         8          &  6864 &  13728 &  41184 &  137280 &  480480 & 1729728 &         &         &         &  \\
    	         9          & 25740 &  51480 & 154440 &  514800 & 1801800 &         &         &         &         &  \\
    	        10          & 97240 & 194480 & 583440 & 1944800 &         &         &         &         &         &  \\ \bottomrule
    \end{tabular}
    \caption{
        The AP bound of graphs consisting of two cliques of size $c_1$ and $c_2$ 
        joint by a single edge.}
    \label{tab:2cluster}
\end{table}

%\begin{equation*}
%    2 \bar{\mu}_i \bar{\mu}_j = 
%    2 \cdot \binom{2(i-1)}{i-1} \cdot \binom{2(j-1)}{j-1}
%\end{equation*}

%\FloatBarrier
%%===============================================================================
%\subsubsection{Multiple cliques of equal size connected via edges}

Table \ref{tab:mcluster} shows the AP bounds of the more general cases where
the networks consist of multiple completely connected subnetworks of the same 
sizes connected via edges to form chains-like structure.
See, for example, the network shown in Figure \ref{fig:mcluster} where five
cliques each of size three are connected via edges that, on a macro level,
resembles a chain.

\begin{table}[h!]
    \centering
    \tiny
    \begin{tabular}{crrrrrrrr}
    	\toprule
    	$c \backslash m$ &    1 &       2 &       3 &       4 &      5 &       6 &        7 &     8 \\ \midrule
    	                 &    1 &       2 &       3 &       4 &      5 &       6 &        7 &     8 \\
    	       1         &      &       2 &       4 &       8 &     16 &      32 &       64 &   128 \\
    	       2         &    2 &       8 &      32 &     128 &    512 &    2048 &     8192 & 32768 \\
    	       3         &    6 &      72 &     864 &   10368 & 124416 & 1492992 & 17915904 &  \\
    	       4         &   20 &     800 &   32000 & 1280000 &        &         &          &  \\
    	       5         &   70 &    9800 & 1372000 &         &        &         &          &  \\
    	       6         &  252 &  127008 &         &         &        &         &          &  \\
    	       7         &  924 & 1707552 &         &         &        &         &          &  \\
    	       8         & 3432 &         &         &         &        &         &          &  \\ \bottomrule
    \end{tabular}
    \caption{
        The AP bound for graphs consisting of $n$ cliques 
        each of size $c$.
        }
    \label{tab:mcluster}\vspace{-0.3in}
\end{table}

\FloatBarrier

\section{IEEE 14 bus system}\label{sec:ieee14}
\begin{wrapfigure}{h}{0.40\textwidth}\vspace{-0.2in}
    \centering
    \includegraphics[width=0.39\textwidth]{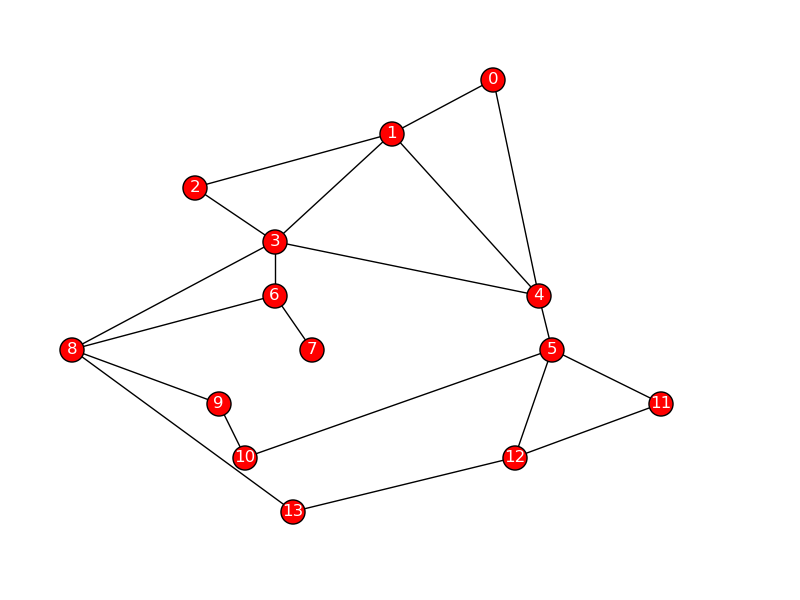}
    \caption{IEEE 14 bus system}
    \label{fig:ieee14}\vspace{-0.3in}
\end{wrapfigure}
The IEEE 14 bus system (the network toopology shown in Figure \ref{fig:ieee14}), representing a portion of the power system of the 
Midwestern U.S.A. in the 1960s, is a widely used benchmark system in testing
the efficiency of different solvers for load flow equations.

Here we show that the isolated complex solution count, BKK bound, and AP bound
are much smaller than previously studied solution bounds.
In particular, the BKK bound in our formulation of the load flow equations is 427680. This means the polyhedral homotopy
method described in Sec.~\ref{sec:sol_count_homotopy} need to trace at most
427680 paths to obtain \textit{all} isolated non-deficient complex solutions, which compared with the BKK bound in 
a previous formulation of the load flow equations \cite{mehta2014a}, 49283072, is around a factor 115 of reduction. 
This is well within the reach of modern polyhedral homotopy implementations.
In particular, with a random choice of the $Y$-matrix, \tech{Hom4PS-3} 
\cite{chen_hom4ps3,chen_hom4ps3_2014} was able to find all solutions in less
than 5 minutes (around 297 seconds) on a single workstation 
(with two \tech{Intel Xeon} processor).
\begin{wraptable}{l}{0.45\textwidth}\vspace{-0.3in}
    \centering
    \begin{tabular}{lr}
        \toprule
        Solutions     & 427680 \\
        BKK bound     & 427680 \\
        AP bound      & 427680 \\
        BKK-MNT\cite{mehta2014a} & 49283072\\
        BBLSY bound   & 10400600 \\
        CB bound      & 67108864 \\ \bottomrule
    \end{tabular}
    \caption{Solution bounds for the IEEE 14-bus topology, with generic complex coefficients. BKK-MNT is the BKK bound with a different 
    polynomial formulation \cite{mehta2014a}.}
    \label{tab:ieee14}\vspace{-1in}
\end{wraptable}

%\noindent

%===============================================================================
\section{Conclusion and Outlook}\label{sec:conclusion}
Solutions of the load flow equations play an important role in design and operations of power systems
and a huge amount of research activities have been focused on developing efficient computational 
methods to solve the load flow equations \cite{mehta2015recent}. 
In this paper we focus on 
computing a tight upper bound on the number of load flow solutions for the following reasons:
the NPHC method \cite{sommese_numerical_2005,li_numerical_2003} introduced in power systems areas in 
\cite{mehta2014a,mehta2014b,chandra2015equilibria,mehta2015algebraic,molzahn-mehta-matt2015} 
requires an upper bound on the number of isolated (complex) solutions of the system. The CB
bound or
most of the other existing upper bounds, and in particular
the BBLSY bound \cite{baillieul1982,li_numerical_1987}, do not capture the 
network topology\cite{molzahn-mehta-matt2015} of the power system into account.
To make the computation of the NPHC method efficient, and to provide a definitive stopping criterion 
for other iterative methods, one has to come up with a tighter upper bound on the number of isolated 
complex solutions.

In this paper, for technical reasons, we ignore the physically less likely 
solutions at which at least one of the voltages is zero. We described a specific
version of the algebraic formulation of the load flow equations and a corresponding tighter upper bound known as the BKK bound.
We showed that there exists at least some generic parameter values for which the BKK bound is attainable. 

Another novel feature in this paper is the introduction of a novel bound, called an AP, 
for sparse systems. With this bound, we proved one of 
our main theorems, Theorem \ref{thm:polytope-bound}, which gives a more explicit bound on the number of complex solutions.
The new bound is computationally much more tractable for large systems compared to the BKK bound.

Our bounds not only give the precise number of complex solutions
for the given power system, but also provide a particular specialization of the NPHC method that in turn
achieves all the solutions for the
system but more efficiently than the previous ones that are based on the crude bounds.
We examined our results on upper bounds on the number of complex solutions with the actual number of complex solutions obtained with the help of 
the polyhedral NPHC method for various cases including complete graphs; path graphs; ring graphs; wheel graphs; random tree graphs; graphs 
having two cliques sharing exactly one, two and three nodes; graphs having two cliques both connected by a single edge; graphs with multiple 
cliques of equal size connected via edges, and the IEEE standard 14-bus system. In all the instances we have tested, 
we find that (1) the BKK bound is always equal to the adjacency bound,
and (2) the BKK bound is always equal to the actual number of complex solutions, i.e., the BKK 
or the AP bound are the actual number of isolated complex solutions for generic coefficients
and completely capture the specific network topology of the given power systems. A proof supporting the numerical results
will be published elsewhere.
Our numerical investigations reproduced all the previously known results \cite{baillieul1982,guo1994, molzahn-mehta-matt2015}. 
In addition, we provide further novel results for aforementioned list of graphs.

We emphasis that the BKK bound
is dependent on the polynomial formulation used to compute it, though the number of complex solutions and number 
of physical solutions are obviously independent of the polynomial formulations. The phrase "BKK bound" throughout the paper, unless 
specified otherwise, is computed from the polynomial formulation we have described in this paper (given in eq. (\ref{equ:powerflow-alg})). 
In other words, 
our polynomial formulation seems optimal since the BKK bound computed using it is generically equal to the number of isolated complex solutions.
To give an example, for the IEEE standard 14-bus system, the BKK using a previous formulation was $49,283,072$, whereas the BKK in 
the present formulation is $427,680$, i.e., a factor of around $115$ reduction with our present formulation.

%\todo[inline]{
%    @Dhagash: I have removed the last three paragraphs of this section.
%    Those are the topics that I know very little about and hence feel quite
%    uncomfortable with:
%    The complexity issue; The future directions; and the ``tracing-method''.
%    }
%Another major implication of our tight bound is that with the aid of it the NPHC method becomes highly efficient in 
%solving the load flow equations, at least for the generic load flow parameters:
%for the NPHC method, each converging and nonsingular path is shown to have ~$O(N^3)$ complexity\cite{sommese_numerical_2005}.
%Since our numerical investigations yield that the BKK bound is sharp in this case, the complexity 
%of each of the paths (the total number of which is equal to 
%the BKK bound of the system) is ~$O(N^3)$ making it computationally efficient to find all the solutions for large networks.

Our results open up several future directions: 
though a less likely case, extending our results to the situtation 
where at least one of the voltages is zero would further complete the discussion on number of complex solutions.
Further refinement will be accomplished by computing or estimating the number of physical 
solutions, and analysing their stability properties for different network topologies.
Investigating for which parameter 
values the number of physical solutions is equal or at least of the same order of magnitude as the 
number of complex solutions of the systems will also provide a better tool to analyse efficiencies of 
different computational methods: For example, for the 3 nodes complete graph 
case, it is known that the number of real solutions does attain the tight-most upper bound on the number of 
complex solutions, $6$\cite{tavora1972equilibrium,baillieul1984critical,klos1991physical}. 
However, for $4$ nodes complete graph case, so far parameter points with only $14$ physical solutions are known as 
opposed to $20$ complex solutions \cite{baillieul1982}.
This has an important consequences on the tracing-methods for solving load flow equations such as proposed 
in Refs.~\cite{thorp1993,Lesieutre2015}. Leaving aside the lack of 
proof for guaranteeing all the physical solutions of the load flow equations in the tracing-methods,
for such parameter points at which the number of physical solutions is equal or of the same order of magnitude as the BKK bound,
the number of traces the method has to go through will be equivalent to the polyhedral NPHC method will also have to go through.
Hence, in such cases, there will be no gain even if one ignores the guarantee to find all solutions that the NPHC method provides.
On the other hand, our bounds can provide the necessary stopping criterion to the tracing-methods and hence the two methods 
may complement each other in some cases.

\section*{Acknowledgement}
TC was supported by NSF under Grant DMS 11-15587. DBM was support by NSF-ECCS award ID 1509036 and an Australian Research
Council DECRA fellowship no. DE140100867. 
%===============================================================================
\appendix

%===============================================================================
\section{Mixed volume and BKK bound} \label{sec:mvol}
Given $m$ convex polytopes $Q_1,\dots,Q_m \subset \R^n$ and $m$ positive 
real numbers $\lambda_1,\dots,\lambda_m$ with $m \le n$ Minkowski's Theorem 
\cite{minkowski_theorie_1911} states that the $n$-dimensional volume of the 
\emph{Minkowski sum} \cite{minkowski_theorie_1911,groemer_minkowski_1977}.
\[
    \lambda_1 Q_1 + \cdots + \lambda_n Q_n =
    \{ \;
        \lambda_1 \boldq_1 + \cdots + \lambda_n \boldq_n \mid
        \boldq_i \in Q_i \text{ for } i = 1,\dots,n \;
    \}
\]
is a homogeneous polynomial of degree $n$ in the variables 
$\lambda_1,\dots,\lambda_n$.
For a $m$-tuple of nonnegative integers $(k_1,\dots,k_m)$, the coefficient 
associated with the monomial $\lambda_1^{k_1} \cdots \lambda_m^{k_m}$ is known 
as the \emph{mixed volume} \cite{groemer_minkowski_1977,minkowski_theorie_1911} 
of type $(k_1,\dots,k_m)$.
Mixed volume of type $(1,\dots,1)$ is often simply referred to as the mixed volume.
As a common convention, the mixed volume involving two shapes in $\R^2$ are
usually called \emph{mixed area} instead.

\begin{figure}[h]
    \centering
    \begin{subfigure}[b]{0.2\columnwidth}
        \centering
        \begin{tikzpicture}[scale=1.0]
        \draw [<->] (1.5,0) -- (0,0) -- (0,1.5);
        \draw [thick] 
        (0,0) node[below] {$a$} -- 
        (1,0) node[below] {$b$} -- 
        (1,1) node[right] {$c$} -- 
        (0,1) node[left ] {$d$} -- 
        cycle;
        \end{tikzpicture}
        \caption{$Q_1$}
        \label{fig:newt-example-a}
    \end{subfigure}%
    ~ 
    \begin{subfigure}[b]{0.22\columnwidth}
        \centering
        \begin{tikzpicture}[scale=1.0]
        \draw [<->] (2.5,0) -- (0,0) -- (0,2.5);
        \draw [thin,dashed] (0,0) grid (2,2);
        \draw [thick] 
        (0,0) node[below] {$e$} -- 
        (2,1) node[right] {$f$} -- 
        (1,2) node[above] {$g$} -- 
        cycle;
        \end{tikzpicture}
        \caption{$Q_2$}
        \label{fig:newt-example-b}
    \end{subfigure}%
    \begin{subfigure}[b]{0.57\columnwidth}
        \centering
        \begin{tikzpicture}[scale=1.1]
        \draw [<->] (3.5,0) -- (0,0) -- (0,3.5);
        \draw [pattern=north west lines, pattern color=blue] 
        (0,0) node[below] {$a+e$} -- 
        (1,0) node[below] {$b+e$} -- 
        (1,1) -- 
        (0,1) node[left ] {$d+e$} -- 
        cycle;
        \draw [pattern=dots, pattern color=red]
        (1,1) -- 
        (3,2) node[right] {$c+f$} -- 
        (2,3) node[above] {$c+g$} -- 
        cycle;
        \draw [fill=lightgray]
        (0,1) --
        (1,3) node[above] {$d+g$} -- 
        (2,3) --
        (1,1) --
        cycle;
        \draw [fill=lightgray]
        (1,0) --
        (3,1) node[right] {$b+f$} -- 
        (3,2) --
        (1,1) --
        cycle;
        \draw [thin,dashed] (0,0) grid (3,3);
        \end{tikzpicture}
        \caption{The Minkowski sum $Q_1 + Q_2$}
        \label{fig:newt-example-c}
    \end{subfigure}
    \caption{
        The Newton polygons of $p_1$ and $p_2$ along with their Minkowski sum
        and a graphical depiction of the mixed volume (mixed area).}
    \label{fig:newt-example}
\end{figure}
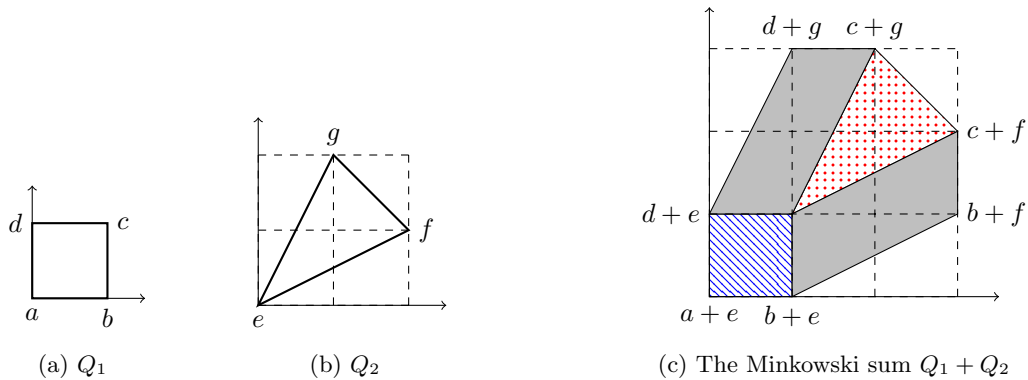

%\section{Newton polytopes and the BKK bound} \label{sec:newt}
%In this section we give a brief explaination of the concepts of Newton polytope
%and its relation to the BKK bound used in Theorem \ref{thm:bernshtein}.
The BKK bound is formulated in terms of mixed volume.
We start with an example (used in \cite{huber_polyhedral_1995}):
With $\boldx=(x_1, x_2)$, let $P(\boldx) $ be the polynomial system
\begin{align*}
    p_1(\boldx)  &= c_{11}x_1x_2   + c_{12}x_1      + c_{13} x_2 + c_{14} \\
    p_2(\boldx)  &= c_{21}x_1x_2^2 + c_{22}x_1^2x_2 + c_{23} 
\end{align*}
Here, $c_{ij}$'s are nonzero coefficients. 
The formal expressions for the monomials in $p_1$ are $x_1x_2 = x_1^1x_2^1$, 
$x_1 = x_1^1 x_2^0$, $x_2 = x_1^0x_2^1$ and $1=x_1^0x_2^0$. 
Considering the exponents of each monomial as a point on the plane, 
the set of exponent points $S_1=\{a=(0,0),\; b=(1,0), \; c=(1,1),\; d=(0,1) \}$
is called the \emph{support} of $p_1$, and its convex hull $Q_1=\op{conv}(S_1)$ 
is called the \emph{Newton polygon} of $p_1$. 
Similarly, the support of $p_2$ is $S_2=\{e=(0,0),\; f=(2,1),\; g=(1,2)\}$ 
and its Newton polygon is $Q_2=\op{conv}(S_2)$.
The two newton polygons are shown in Figure \ref{fig:newt-example-a} and
\ref{fig:newt-example-b} respectively.
The mixed area of $Q_1$ and $Q_2$ can be intuitively taken as the area within
the Minkowski sum $Q_1 + Q_2$ contributed by ``mixing'' of the two shapes.
From Figure \ref{fig:newt-example-c}, it is easy to see that $Q_1 + Q_2$
contains copies of $Q_1$ and $Q_2$. The rest (the shaded region) are resulted
from the ``mixing effect'' of edges from both $Q_1$ and $Q_2$ in constructing
the Minkowski sum. The total area of these regions is the mixed area.

%With the notation $\boldx^{\boldq} = x_1^{q_1} x_2^{q_2}$ for $\boldq=(q_1,q_2)$,
%we may rewrite \eqref{AA} as
%\[
%p_1(\boldx) = \sum_{\boldq \in S_1} c_{1,\boldq} \boldx^{\boldq} \qquad \text{and} \qquad
%p_2(\boldx) = \sum_{\boldq \in S_2} c_{2,\boldq} \boldx^{\boldq}.
%\]

In general, given a polynomial system $P = (p_1,\dots,p_n)$ of the form
\begin{equation*}
    \left(
        \sum_{\bolda\in S_1}c_{1,\bolda}\boldx^{\bolda},
        \;\dots\;,
        \sum_{\bolda\in S_n}c_{n,\bolda}\boldx^{\bolda}
    \right)
\end{equation*}
where $\boldx=(x_1,\dots,x_n)$ and $S_1,\dots,S_n$ are fixed subsets of $\Z^n$ 
called the {\em support} of $p_1,\dots,p_n$ respectively. 
Note that here with $\bolda=(a_1,\dots,a_n)$ we use the notation
$\boldx^{\bolda} =x_1^{a_1} \cdots x_n^{a_n}$ for a monomial.
The convex hull $Q_j=\op{conv}(S_j)$ in $\R^n$ is called the {\em Newton polytope} 
of $p_j$. 
In this context, the BKK bound is the mixed volume of the Newton polytopes
$Q_1,\dots,Q_n$.

%===============================================================================
\section{Zariski openness} \label{sec:zariski}
Part (B) of Theorem \ref{thm:bernshtein} states that ignoring the conjugate relations among 
the coefficients, the set of coefficients for which the upper bound given by
Theorem \ref{thm:bernshtein} part (A) hold exactly is at least a nonempty 
``Zariski open'' set.
This ``Zariski openness'' can be characterized as follows:
If we replace $Y_{jk}^*$ and $S_j^*$ in \eqref{equ:powerflow-alg} by independent 
coefficients $Y_{jk}'$ and $S_j'$, then there exists a nonzero polynomial $D$ in 
the variables $\{ Y_{jk},Y_{jk}',S_j,S_j' \}$
such that $D \ne 0$ implies the BKK bound is exact.
    
%===============================================================================
\section{The polarization lemma} \label{sec:polarization}
The proof for Theorem \ref{thm:attainable} hinges on a key lemma in the theory
of complex variables that is sometimes referred to as the 
\textit{polarization lemma}:

\smallskip
\begin{lemma}
    \label{lem:polarization}
    Suppose that $H : \Cstar^n \times \Cstar^n \to \C$ is a holomorphic function
    of the $2n$ complex variables $(\boldz,\mathbf{w})$, and that $H(\boldz,\boldz^*) = 0$
    for all $\boldz \in \Cstar^n$. Then, $H(\boldz,\mathbf{w}) = 0$
    for all $(\boldz,\mathbf{w}) \in \Cstar^n \times \Cstar^n$.
\end{lemma}
\smallskip

This lemma has its roots in the works of W. Wirtinger.
Here we refer to \cite{dangelo_several_1993} for a direct proof.

\bibliographystyle{siam}
\bibliography{power-flow,acc}
\end{document}